\documentclass[reqno]{amsart}

\usepackage{amssymb}
\usepackage{amsfonts}
\usepackage{mathrsfs}

\usepackage{xcolor}

\theoremstyle{plain}
    \newtheorem{theorem}{Theorem}[section]
    \newtheorem*{theorem*}{Theorem}
    \newtheorem{lemma}[theorem]{Lemma}
    \newtheorem{proposition}[theorem]{Proposition}

\theoremstyle{definition}
    \newtheorem{definition}{Definition}[section]
    \newtheorem{remark}{Remark}[section]
    
\numberwithin{equation}{section}

\renewcommand{\r}{\right}
\newcommand{\m}{\mbox{\textnormal{mid}}}
\begin{document}

\title[]{Nonexistence of $(1,6)$-sign soliton solutions for the damped nonlinear Klein--Gordon equation in 2D}
\author[K. Ishizuka]{Kenjiro Ishizuka
}
\address[K. Ishizuka]{Azabu Junior and Senior High School, Tokyo, Japan}
\email{k-ishizuka@azabu-jh.net}
\date{}
\keywords{damped nonlinear Klein-Gordon equation, multi-solitons, long-time asymptotics, interaction dynamics, asymptotic geometry, rigidity}
\subjclass[2020]{35L71, 35B40, 35B44}

\begin{abstract}
We consider the two-dimensional damped nonlinear Klein--Gordon equation
\begin{align*}
\partial_t^2u-\Delta u+2\alpha\partial_tu+u-|u|^{p-1}u=0
\end{align*}
on $\mathbb{R}^2$, where $\alpha>0$ and $p>2$.
We study multi-soliton solutions asymptotic in the energy space, in the full-time limit $t\to\infty$, to a superposition of translated ground states, with one soliton having the sign opposite to that of the other $k$ solitons.
C\^{o}te and Du conjectured nonexistence for $k\geq7$ in the regular polygonal setting and identified the hexagonal case $k=6$ as critical.
In this paper, we prove that $(1,k)$-sign soliton solutions do not exist for any $k\geq6$, without imposing any symmetry assumptions.
\end{abstract}

\maketitle

\tableofcontents

\section{Introduction}
\subsection{Setting of the problem}
We consider the following damped nonlinear Klein--Gordon equation
\begin{align}
\label{DNKG}
\tag{DNKG}
	\left\{
	\begin{aligned}
		&\partial_{t}^2u-\Delta u+2\alpha \partial_{t}u+u-f(u)=0, & (t,x) \in \mathbb{R} \times \mathbb{R}^2,
		\\
		&\left(u(0,x),{\partial}_tu(0,x)\right)=\left(u_{0}(x), u_1(x)\right)\in \mathcal{H},
	\end{aligned}
	\r.
\end{align}
where $\alpha>0$, $\mathcal{H}=H^1(\mathbb{R}^2)\times L^2(\mathbb{R}^2)$, $f(u)=|u|^{p-1}u$, with a power $p>2$. It follows from \cite{BRS} that the Cauchy problem for \eqref{DNKG} is locally well-posed in the energy space $\mathcal{H}$. Moreover, defining the energy of a solution $\vec{u}=(u, {\partial}_tu)$ by
\begin{align*}
E(\vec{u}(t))=\frac{1}{2}\|\vec{u}(t)\|_{\mathcal{H}}^2-\frac{1}{p+1}\| u(t)\|_{L^{p+1}}^{p+1},
\end{align*}
we have
\begin{align}\label{energydecay}
E(\vec{u}(t_2))-E(\vec{u}(t_1))=-2\alpha \int_{t_1}^{t_2} \| {\partial}_tu(t)\|_{L^2}^2dt.
\end{align}
In addition, the equation \eqref{DNKG} enjoys several invariances:
\begin{align}\label{invariance}
\begin{aligned}
u(t,x)&\mapsto -u(t,x),
\\
u(t,x)&\mapsto u(t,x+y),
\\
u(t,x)&\mapsto u(t+s,x),
\\
u(t,x)&\mapsto u(t,Rx),
\end{aligned}
\end{align}
where $y\in \mathbb{R}^2$ and $s\in \mathbb{R}$, and 
\begin{align*}
R\in SO(2):=\left\{A\in\mathbb{R}^{2\times 2}:A^{\mathsf T}A=I,\ \det A=1\right\}.
\end{align*}

\subsection{Long-time asymptotics of multi-soliton solutions}

In this paper, we are interested in the dynamics of multi-soliton solutions related to the ground state $Q$, which is the unique positive radial solution in $H^1(\mathbb{R}^2)$ of
\begin{align}\label{Qode}
-\Delta Q+Q-f(Q)=0.
\end{align}
See \cite{K2}. The ground state generates the stationary solution $(Q,0)$. By the invariances in \eqref{invariance}, $(-Q,0)$ and every spatial translate $(Q(\cdot+y),0)$ are also stationary solutions of \eqref{DNKG}.

There has been intensive study of the global behavior of general, possibly large, solutions to nonlinear dispersive equations. A guiding principle is the \textit{soliton resolution conjecture}, which predicts that generic global solutions are asymptotic to a superposition of solitons for large time. For \eqref{DNKG}, Feireisl \cite{F} established such a decomposition along a sequence of times for every bounded global solution. In one space dimension, where all nonzero stationary solutions are translates of $\pm Q$, C\^{o}te, Martel, and Yuan \cite{CMY} proved soliton resolution in the full-time limit and showed that the solitons must have alternating signs.

Even when a soliton decomposition is known, the decomposition alone does not determine which collections of solitons can actually arise or what constraints the dynamics imposes on their parameters. This naturally leads to the problem of classifying the admissible asymptotic configurations of the constituent solitons. In this paper, we consider superpositions of translates of $\pm Q$, so the relevant parameters are the signs and the time-dependent soliton centers.

Understanding the dynamics of these centers requires controlling both the mutual interactions among the solitons and the radiative remainder; see, for example, \cite{CJ,IN}. However, for \eqref{DNKG}, the energy dissipation and the decay
of the linear damped flow provide control of the radiative component, while the damping causes the soliton velocities
to decay. Consequently, the long-time evolution of the centers is governed, to leading order, by the mutual soliton interactions and is described by an effective first-order system.

This damping-induced reduction has enabled considerable progress in this direction for \eqref{DNKG}. C\^{o}te, Martel, Yuan, and Zhao \cite{CMYZ} classified the case of two solitons in general space dimension: the two solitons must have opposite signs, and their centers separate along a fixed straight line. Their analysis also makes explicit the leading-order interaction law, according to which like-signed solitons attract and oppositely signed solitons repel.  Combining the sign restrictions obtained in \cite{CD} with the results of the author \cite{I2}, the case of three solitons is also classified: exactly one soliton has the sign opposite to that of the other two, and the three centers are necessarily asymptotically collinear. More recently, the author \cite{I3,I4} has classified the long-time asymptotics of four-soliton solutions associated with the ground state. In \cite{I3}, the author proved that, when one soliton has the sign opposite to that of the other three, its center converges to a fixed point, while the other three centers form an expanding asymptotically equilateral triangle around it. In \cite{I4}, the author proved that, when there are two positive and two negative solitons, the centers asymptotically form either an expanding collinear configuration or an expanding rhombus, with alternating signs in both cases.

In this paper, we focus on multi-soliton solutions in two space dimensions that are asymptotic to a superposition of ground states with exactly one soliton having the sign opposite to that of all the others. By imposing suitable symmetry assumptions, C\^{o}te and Du \cite{CD} constructed solutions in which an oppositely signed soliton remains at the center, while the other soliton centers form an expanding equilateral triangle, square, or regular pentagon. In these configurations, the attractive interactions among the outer solitons are negligible compared with their repulsive interactions with the central soliton, which makes the construction possible. For a regular hexagon or heptagon, however, these attractive interactions can no longer be neglected, making the construction more delicate.

Indeed, in \cite[Remark 1.7]{CD}, C\^{o}te and Du conjectured the nonexistence of such multi-solitons for regular polygons with $k\geq7$ vertices and identified the hexagonal case $k=6$ as critical. In this paper, we prove, without imposing any symmetry assumptions, that no global solution is asymptotic in the energy space, in the full-time limit $t\to\infty$, to a superposition of seven or more mutually separating ground states with exactly one soliton having the sign opposite to that of all the others.

\subsection{Main result}
First, we define solutions that behave as the superposition of ground states.
\begin{definition}\label{defKsoli}
A global solution $\vec{u}$ of \eqref{DNKG} is called a \textit{$K$-soliton solution} for $K\in\mathbb{N}$ if there exist $\sigma \in \{-1,1\}^K$ and $z:[0,\infty)\to (\mathbb{R}^2)^K$ such that 
\begin{align}
\label{Ksoli}
\begin{aligned}
\lim_{t\to \infty} \left(\|u(t)-\sum_{i=1}^K\sigma_i Q(\cdot-z_i(t))\|_{H^1}+\|{\partial}_tu(t)\|_{L^2}\right)&=0,
\\
\lim_{t\to\infty}\left( \min_{i\neq j}|z_i(t)-z_j(t)|\right)&=\infty.
\end{aligned}
\end{align}
\end{definition}
\begin{remark}
When $K=1$, we adopt the convention that $\min_{i\neq j}|z_i(t)-z_j(t)|=\infty$.
\end{remark}
\begin{remark}
In \cite[Definition 1.11]{CD}, a $K$-soliton solution is defined in terms of convergence along a sequence of times $t_n\to\infty$, whereas Definition~\ref{defKsoli} requires convergence in the full-time limit $t\to\infty$.
Our nonexistence theorem is proved under this full-time assumption; nonexistence under the subsequential assumption is not established here.
\end{remark}

In this paper, we focus on $K$-soliton solutions in which exactly one soliton has the sign opposite to that of all the others. To specify the numbers of copies of $Q$ and $-Q$ in the soliton decomposition, we introduce the notion of $(m,n)$-sign soliton solutions as follows.
\begin{definition}\label{defm,nsoli}
A global solution $\vec{u}$ of \eqref{DNKG} is called a \textit{$(m,n)$-sign soliton solution} for $m,n\in\mathbb{N}\cup\{0\}$ if there exist $z_+:[0,\infty)\to (\mathbb{R}^2)^m$ and $z_-:[0,\infty)\to(\mathbb{R}^2)^n$ such that 
\begin{align}
\label{mnsoli}
\begin{aligned}
\lim_{t\to \infty} \left(\|u(t)-\sum_{i=1}^mQ(\cdot-z_{+,i}(t))+\sum_{i=1}^n Q(\cdot-z_{-,i}(t))\|_{H^1}+\|{\partial}_tu(t)\|_{L^2}\right)&=0,
\\
\lim_{t\to\infty}\left( \min_{(\circ,i)\neq (\bullet,j)}|z_{\circ,i}(t)-z_{\bullet,j}(t)|\right)&=\infty.
\end{aligned}
\end{align}
Here, the minimum is taken over all distinct elements $(\circ,i)$ and $(\bullet,j)$ of $\{(+,i):1\leq i\leq m\}\cup\{(-,i):1\leq i\leq n\}$.
\end{definition}
\begin{remark}
As in Definition \ref{defKsoli}, when $m+n\leq1$, we adopt the convention that $\min_{(\circ,i)\neq(\bullet,j)}|z_{\circ,i}(t)-z_{\bullet,j}(t)|=\infty$.
\end{remark}
\begin{remark}
Since \eqref{DNKG} is invariant under the sign change $u\mapsto -u$, any property established for $(m,n)$-sign soliton  solutions immediately yields the corresponding property for $(n,m)$-sign soliton solutions.
\end{remark}
As mentioned above, we focus on $(1,k)$-sign soliton solutions with $k\geq6$. Our main theorem states that no such solutions exist.
\begin{theorem}\label{maintheorem}
For $k\geq 6$, there are no $(1,k)$-sign or $(k,1)$-sign soliton solutions.
\end{theorem}

\subsection{Main difficulties and new ideas}

Determining the asymptotic configuration of multi-soliton solutions with many components is delicate because the attractive and repulsive interactions can compete and cancel. In \cite[Section 4]{I4}, under a suitable repulsive condition, the author used an argument based on the deficit in the Cauchy--Schwarz inequality to establish subsequential convergence of the suitably normalized soliton centers to a configuration satisfying balance equations. However, this repulsive condition is not readily derived for $(1,k)$-sign soliton solutions with $k\geq6$. In particular, when the centers of six solitons of the same sign lie at the vertices of a regular hexagon, with an oppositely signed soliton at its center, the repulsive force exerted by the central soliton on each outer soliton is canceled at leading order by the combined attractive forces exerted by its two adjacent outer solitons. 

To overcome this difficulty, we first refine the analysis of the soliton interactions, retaining the higher-order pair interactions and obtaining sharper estimates for the remainder involving three or more solitons. We then establish quantitative estimates for perturbations of a cluster of six soliton centers from a regular hexagon around the oppositely signed soliton. These estimates relate the interaction energy to the forces generated by deviations from the hexagonal configuration.

Combining these perturbation estimates with the energy estimates and the dissipative dynamics, we derive a positivity estimate for the interaction energy corresponding to the repulsive condition. This allows us to apply the deficit argument of \cite{I4} and obtain the limiting balance equations. The final contradiction follows from an elementary argument in planar geometry. 

\subsection{Previous results}

The present work lies at the intersection of the soliton resolution conjecture and the study of strongly interacting multi-soliton dynamics. For the focusing energy-critical wave equation, soliton resolution has been established along sequences of times for general bounded solutions in dimensions three to five; see \cite{DJKM}. Full-time soliton resolution for radial solutions with bounded energy norm is now known in a broad range of dimensions; see, in particular, \cite{DKM2,CDKM,JL4}. Related full-time results have also been obtained for equivariant wave maps \cite{DKMM,JL3}, for radial energy-critical wave equations with a potential \cite{JLX,LMZ}, for radial energy-critical wave equations with damping \cite{GZ}, and for the radial energy-critical heat equation \cite{A1}.

More recently, the rigidity of geometric parameters within multi-bubble regimes has also been investigated. Kim and Merle \cite{KM2} classified several non-radial multi-bubble regimes for the energy-critical nonlinear heat equation in dimensions $N\geq7$, in which the scales, centers, and signs of the bubbles enter the effective dynamics in a nontrivial manner. For the five-dimensional focusing energy-critical wave equation, Jendrej, Zhang, and Zhao \cite{JZZ} proved that, for solutions asymptotic in the energy space to a sum of positive ground states whose scales remain comparable and whose centers
converge to distinct points, all scaling parameters are of order $t^{-2}$ and the renormalized modulation parameters converge to a finite-dimensional algebraic set determined by the limiting centers.

For the nonlinear Klein--Gordon equation without damping, the soliton resolution conjecture remains widely open.
Nevertheless, Nakanishi and Schlag \cite{NS} obtained a complete classification of the global dynamics for the three-dimensional focusing cubic equation when the energy is at most slightly above that of the ground state. For the same equation, Chen and Jendrej \cite{CJ} established conditional asymptotic stability of well-separated multi-solitons with pairwise distinct velocities and classified the corresponding pure multi-soliton solutions. In the damped setting, Keller \cite{K1} constructed stable and unstable manifolds near stationary solutions. Subsequential decomposition results for bounded global solutions were obtained by Feireisl \cite{F} and Li and Zhao \cite{LZ}, under suitable assumptions on the dimension and the nonlinearity. Full-time results include convergence for radial solutions \cite{BRS} and soliton resolution in one space dimension without any size restriction \cite{CMY}.

For the damped nonlinear Klein--Gordon equation, C\^{o}te, Martel, Yuan, and Zhao \cite{CMYZ} determined the universal logarithmic separation law for two-soliton solutions associated with the ground state and constructed the corresponding family of solutions as a codimension-two Lipschitz manifold in the energy space. Subsequently, the author and Nakanishi \cite{IN} classified the global dynamics of all solutions starting sufficiently close to a superposition of two well-separated ground states with opposite signs. Results concerning solutions associated with excited stationary states can be found in \cite{CY}.

Multi-soliton solutions have also been constructed for several nonlinear dispersive equations. For the nonlinear Klein--Gordon equation without damping, C\^{o}te and Mu\~{n}oz \cite{CM} constructed multi-solitons based on the ground state, and C\^{o}te and Martel \cite{CM1} extended this construction to more general stationary profiles. These solutions are asymptotic to solitary waves traveling with pairwise distinct velocities of magnitude strictly less than one. Strongly interacting multi-solitons with logarithmic relative distances have been obtained for several models; see, for example, \cite{A,GI,MN,N1,N2}. In the damped setting considered here, the soliton velocities decay because of the damping, and the logarithmic separation is instead generated by the balance between the damping and the exponentially small interactions.

\subsection{Notation}

Let $\{e_1,e_2\}$ denote the canonical basis of $\mathbb{R}^2$, and let ${\partial}_k$ denote the partial derivative
with respect to $x_k$. For $x,y\in\mathbb{R}^2$, we denote the Euclidean inner product by
\begin{align}\label{Euclideanproduct}
x\cdot y=x_1y_1+x_2y_2.
\end{align}
For any pairs of functions $u,v\in L^2(\mathbb{R}^2)$, we denote the $L^2$ scalar product by
\begin{align}\label{L2naisekiteigi}
\langle u,v\rangle=\int_{\mathbb{R}^2}u(x)v(x)dx.
\end{align}
Throughout this paper, the comparison notation is used as follows. For nonnegative quantities $X$ and $Y$, we write $X\lesssim Y$ if there exists a constant $C>0$ such that $X\leq CY$. We write $X\ll Y$ if $X\leq cY$ for a sufficiently small constant $c>0$. The precise smallness requirement will be clear from the context or specified when needed.

For real-valued quantities $X$ and $Y$ with $Y\neq0$, we write $X\sim Y$ if there exist constants $0<c\leq C$ such that
\begin{align}\label{signedcomparison}
c\leq\frac{X}{Y}\leq C.
\end{align}
Thus, $X\sim Y$ means that $X$ and $Y$ have the same sign and comparable absolute values; in particular, this notation is also used when both quantities are negative. For nonnegative $X$ and $Y$, this
is equivalent to $X\lesssim Y$ and $Y\lesssim X$. Unless otherwise stated, implicit constants are independent of the variables under consideration and may change from line to line.

Moreover, in this paper we deal with the second smallest value among three quantities. We therefore define
\begin{align}\label{middef}
\m{(a,b,c)}=a+b+c-\max{(a,b,c)}-\min{(a,b,c)}.
\end{align}

\section{Preliminaries}

In this section, we collect several preliminary results on the long-time dynamics of $K$-soliton solutions. Most of the material is taken from \cite[Section 2 and Section 3]{I2} and is included here for the reader's
convenience.

\subsection{Basic properties of the ground state}

Since $Q$ is radial and positive, there exists $q:[0,\infty)\to (0,\infty)$ such that
\begin{align*}
Q(x)=q(|x|).
\end{align*}
Furthermore, $q$ satisfies
\begin{align}\label{qode2}
q^{\prime \prime}+\frac{1}{r}q^{\prime}-q+q^p=0.
\end{align}
In addition, there exists a constant $c_q>0$ depending only on $p$ such that, for $r\geq 1$,
\begin{align}\label{qes}
|q(r)-c_qr^{-\frac{1}{2}}e^{-r}|+|q^{\prime}(r)+c_qr^{-\frac{1}{2}}e^{-r}|
\lesssim r^{-\frac{3}{2}}e^{-r}.
\end{align}
Next, in order to make the leading interaction term more explicit, we introduce the vector field
$H:\mathbb{R}^2\setminus\{0\}\to\mathbb{R}^2$ defined by
\begin{align}\label{Hdef}
H(z)=-\int_{\mathbb{R}^2}(\nabla Q^p(x))\,Q(x-z)\,dx.
\end{align}
Moreover, we define $g:(0,\infty)\to\mathbb{R}$ as
\begin{align}\label{gdef}
H(z)=\frac{z}{|z|}\,g(|z|).
\end{align}
We note that this definition is well-defined. Furthermore, we have 
\begin{align*}
g(r+h)-g(r)
=-\int_{\mathbb{R}^2} {\partial}_1(Q^p)(x)\Bigl( Q\bigl(x-(r+h)e_1\bigr)-Q(x-re_1)\Bigr)\,dx.
\end{align*}
Therefore, we have 
\begin{align*}
g^{\prime}(r)
=\int_{\mathbb{R}^2} {\partial}_1(Q^p)(x)\,{\partial}_1Q(x-re_1)\,dx
=-\int_{\mathbb{R}^2}Q^p(x)\,{\partial}_1^2Q(x-re_1)\,dx.
\end{align*}
In addition, by direct computation, we have
\begin{align*}
{\partial}_1^2Q(x-re_1)
=\frac{x_2^2}{|x-re_1|^3}q^{\prime}(|x-re_1|)
+\frac{(x_1-r)^2}{|x-re_1|^2}q^{\prime \prime}(|x-re_1|).
\end{align*}
Using \eqref{qode2} and \eqref{qes}, we obtain, for $r\geq1$,
\begin{align}\label{gnozenkin1}
|g(r)-c_gr^{-\frac{1}{2}}e^{-r}|+|g^{\prime}(r)+c_gr^{-\frac{1}{2}}e^{-r}|
\lesssim r^{-\frac{3}{2}}e^{-r},
\end{align}
where $c_g=c_q\int_{\mathbb{R}^2} Q^p(x)e^{-x_1}dx$.

\subsection{Basic properties of the dynamics of $(1,k)$-sign soliton solutions}

We fix $k\geq 6,\ k\in\mathbb N$ and let $\vec{u}$ be a $(1,k)$-sign soliton solution. Then, there exists a $C^1$ function $z:[0,\infty)\to(\mathbb{R}^2)^{k+1}$ such that
\begin{align}\label{1,ksolicond}
\begin{aligned}
\lim_{t\to\infty}(\|u(t)-Q(\cdot-z_0(t))+\sum_{i=1}^kQ(\cdot-z_i(t))\|_{H^1}+\|{\partial}_tu(t)\|_{L^2})=0,
\\
\lim_{t\to\infty}\min_{0\leq i<j\leq k}|z_j(t)-z_i(t)|=\infty.
\end{aligned}
\end{align}
For $z=(z_0,z_1,\ldots,z_k)$, we define $Q_0=Q(\cdot-z_0)$, and for $i=1,2,\ldots,k$, we define $Q_i=-Q(\cdot-z_i)$. Furthermore, we define for $i=0,1,\ldots, k$,
\begin{align*}
\vec{Q}_i=
\begin{pmatrix}
Q_i\\
0
\end{pmatrix},
\qquad
Q_{\sum}=\sum_{i=0}^kQ_i,\qquad
\vec{Q}_{\sum}=
\begin{pmatrix}
Q_{\sum}\\
0
\end{pmatrix}.
\end{align*}
Furthermore, for a solution $\vec{u}$ of \eqref{DNKG} and $z$, we define
\begin{align*}
\vec{\varepsilon}=
\begin{pmatrix}
\varepsilon\\
\eta
\end{pmatrix}
=\vec{u}-\vec{Q}_{\sum}.
\end{align*}
Then $\vec{\varepsilon}$ satisfies the following system:
\begin{align}\label{varepnoode}
\left\{
	\begin{aligned}
		&{\partial}_t\varepsilon=\eta+R_1,
		\\
		&{\partial}_t\eta=\left(\Delta-1+f^{\prime}(Q_{\sum})\right)\varepsilon-2\alpha \eta+R_2+R_3,
	\end{aligned}
	\r.
\end{align}
where
\begin{align*}
R_1&=\sum_{j=0}^k\dot{z}_j\cdot \nabla Q_j,
\\
R_2&=f(Q_{\sum})-\sum_{j=0}^kf(Q_j),
\\
R_3&=f(Q_{\sum}+\varepsilon)-f(Q_{\sum})-f^{\prime}(Q_{\sum})\varepsilon.
\end{align*}

We also introduce the minimal separation
\begin{align*}
D(t)=\min_{i\neq j}|z_i(t)-z_j(t)|.
\end{align*}

First, by suitably choosing the centers as in \cite{I2}, we can strengthen the conditions in \eqref{1,ksolicond} as follows.
\begin{align}\label{1,ksolicondrevise}
\begin{aligned}
\lim_{t\to\infty}\|\vec{\varepsilon}(t)\|_{\mathcal{H}}=0,
\\
\lim_{t\to\infty}D(t)=\infty,
\\
\int_{\mathbb{R}^2} \left(\eta+2\alpha\varepsilon\right){\partial}_lQ(\cdot-z_i(t))=0\ \mbox{for}\ 0\leq i\leq k,\ 1\leq l\leq 2,\ t\gg1.
\end{aligned}
\end{align}

\subsection{Improved estimates for the remainder of $(1,k)$-sign soliton solutions}

In this subsection, we improve the estimate for the size of the remainder $\vec{\varepsilon}$.
\begin{lemma}\label{vareprefinelem}
Let $\vec{u}$ be a $(1,k)$-sign soliton solution satisfying \eqref{1,ksolicondrevise}. Then, we have for $t\gg 1$
\begin{align}\label{epsiloneskai}
\|\vec{\varepsilon}(t)\|_{\mathcal{H}}\lesssim g(D(t)).
\end{align}
\end{lemma}
\begin{proof}
First, by \cite{CGNT}, the operator
\begin{align*}
-\Delta+1-pQ^{p-1}
\end{align*}
on $L^2$ has a unique negative eigenvalue $-\nu_0^2$. Let $\phi$ be an $L^2$-normalized eigenfunction associated with this eigenvalue. Furthermore, we define $\nu^{\pm},\ a_i^+$ ($0\leq i\leq k$) and $a^+$ as
\begin{align*}
\nu^{\pm}&=-\alpha\pm \sqrt{\alpha^2+\nu_0^2},
\\
a_i^+&=\int_{\mathbb{R}^2}\left(\eta-\nu^-\varepsilon\right)\phi(\cdot-z_i)dx,
\\
a^+&=\sqrt{(a_0^+)^2+(a_1^+)^2+\cdots+(a_k^+)^2}.
\end{align*}
By \cite[Section 3]{I2}, we have
\begin{align}
|\dot{z}|+\left|\frac{d}{dt}a^+-\nu^+a^+\right|&\lesssim g(D), \label{prees1}
\\
\lim_{t\to\infty}a^+(t)&=0, \label{prees2}
\\
\lim_{t\to\infty}\frac{\|\vec{\varepsilon}(t)\|_{\mathcal{H}}^2}{g(D(t))}&=0. \label{prees3}
\end{align}
We note that, although \eqref{prees1} may appear to differ from the corresponding estimate in \cite{I2}, it follows immediately from that estimate and \eqref{prees3}. Furthermore, by \cite[Subsection 3.5]{I2}, there exist a $C^1$ function $\mathcal{G}$ and constants $C_1, \tau>0$ such that 
\begin{align}\label{mathcalGpre}
\mathcal{G}+C_1|a^+|^2\sim \|\vec{\varepsilon}\|_{\mathcal{H}}^2,\ \mathcal{G}^{\prime}+\tau \mathcal{G}\lesssim |a^+|^2+g(D)^2.
\end{align}
Here, let $\delta>0$ be sufficiently small. By \eqref{prees1}, we have for $t\gg 1$ $|\dot{D}|\lesssim g(D)$. Therefore, by \eqref{gnozenkin1}, we have for $t\gg 1$
\begin{align*}
\left|\frac{d}{dt}g(D(t))\right|\leq \delta g(D(t)).
\end{align*}
Therefore, for $1\ll t<s$, we obtain
\begin{align}\label{prees4}
e^{-\delta(s-t)}g(D(t))\leq g(D(s))\leq e^{\delta(s-t)}g(D(t)).
\end{align}
In particular, we note that $\delta\ll \nu^+$. We fix $t\gg 1$. Then, by \eqref{prees1}, for every $T>t$, we have
\begin{align*}
\left|a^+(t)-e^{-\nu^+(T-t)}a^+(T)\right|&\leq \int_t^T e^{-\nu^+(s-t)}\left|\frac{d}{ds}a^+(s)-\nu^+a^+(s)\right|ds
\\
&\lesssim \int_t^T e^{-\left(\nu^+-\delta\right)(s-t)}g(D(t))ds
\\
&\lesssim g(D(t)).
\end{align*}
Since we have $\lim_{T\to\infty}e^{-\nu^+(T-t)}a^+(T)=0$, we obtain
\begin{align}\label{prea+es}
a^+(t)\lesssim g(D(t)).
\end{align}
Then, by \eqref{mathcalGpre}, \eqref{prees4}, \eqref{prea+es}, and $\delta\ll \tau$, we have for $t>T\gg 1$
\begin{align*}
\mathcal{G}(t)-e^{-\tau(t-T)}\mathcal{G}(T)&\lesssim \int_T^t e^{-\tau(t-s)}g(D(s))^2ds
\\
&\lesssim \int_T^te^{-(\tau-2\delta)(t-s)}g(D(t))^2ds\lesssim g(D(t))^2.
\end{align*}
Therefore, by \eqref{prees4}, we obtain
\begin{align}\label{mathcalGprees2}
\mathcal{G}(t)&\lesssim e^{-\tau(t-T)}\left|\mathcal{G}(T)\right|+g(D(t))^2\lesssim g(D(t))^2.
\end{align}
By \eqref{mathcalGpre}, \eqref{prea+es} and \eqref{mathcalGprees2}, we obtain \eqref{epsiloneskai}.
\end{proof}

In particular, Lemma \ref{vareprefinelem} yields
\begin{align}\label{R3es}
\|R_3\|_{L^2}\lesssim \|\varepsilon\|_{H^1}^2\lesssim g(D)^2.
\end{align}

\subsection{Improved interaction estimates}

In this subsection, we derive more precise estimates for $R_2$. We define 
\begin{align}\label{Psidef1}
\Psi(r)=\int_r^{\infty} g(s)ds.
\end{align}
Then, we note that by 
\begin{align*}
\lim_{r\to\infty}\int_{\mathbb{R}^2}Q(x)^pQ(x-re_1)dx=0
\end{align*}
and
\begin{align*}
\frac{d}{dr}\int_{\mathbb{R}^2}Q(x)^pQ(x-re_1)dx&=-\int_{\mathbb{R}^2}Q(x)^p{\partial}_1Q(x-re_1)dx
\\
&=\int_{\mathbb{R}^2}{\partial}_1(Q^p)(x)Q(x-re_1)dx
\\
&=-e_1\cdot H(re_1)=-g(r),
\end{align*}
we have
\begin{align}\label{Psidef2}
\Psi(r)=\int_r^{\infty} g(s)ds=\int_{\mathbb{R}^2}Q(x)^pQ(x-re_1)dx.
\end{align}
Furthermore, by direct computation, we have
\begin{align}\label{Psies0922}
\begin{aligned}
&\Psi(r)=g(r)+O\left(\frac{g(r)}{r}\right),\ \Psi^{\prime}(r)=-g(r),
\\
 &\Psi^{\prime \prime}(r)=g(r)+O\left(\frac{g(r)}{r}\right),\ \Psi^{\prime \prime \prime}(r)=O(g(r)).
\end{aligned}
\end{align}
We next estimate the nonlinear interactions. First, Taylor expansion yields
\begin{align}\label{R2taylores}
\|R_2\|_{L^2}+\sum_{i=0}^k\|\left(f^{\prime}(Q_{\sum})-f^{\prime}(Q_i)\right)\nabla Q_i\|_{L^2}\lesssim g(D).
\end{align}
Furthermore, to refine the analysis, we define $\mathcal{N}_s$ for $s=\pm 1$ by
\begin{align}\label{mathcalN+-}
\mathcal{N}_s(r)=\frac{1}{p+1}\int\left(|Q+sQ(\cdot-re_1)|^{p+1}-Q^{p+1}-Q(\cdot-re_1)^{p+1}\right)-2s\Psi(r).
\end{align}
Then, we have the following lemma.
\begin{lemma}\label{mathcalNeslem}
For $s=\pm 1$ and $r\gg 1$, we have
\begin{align}
\mathcal{N}_s(r)>0, \label{mathcalNsei}
\\
\left|\mathcal{N}_s^{\prime}(r)\right|\lesssim \mathcal{N}_s(r). \label{mathcalNbibunes}
\end{align}
Furthermore, for every $1<\gamma<\min{(2,\frac{p+1}{2})}$, we have
\begin{align}\label{mathcalNgensui}
\mathcal{N}_s(r)\lesssim e^{-\gamma r}.
\end{align}
\end{lemma}
\begin{proof}

First, for $s=\pm 1$ and $a, b\geq 0$, we define
\begin{align}
J_s(a,b)=|a+sb|^{p+1}-a^{p+1}-b^{p+1}-s(p+1)(a^pb+ab^p).
\end{align}
Then, by direct computation, we have
\begin{align*}
{\partial}_a{\partial}_bJ_s\geq 0,\ J_s(a,0)=J_s(0,b)=0.
\end{align*}
The inequality is strict when $a,b>0$. Therefore, we have for $a,b>0$
\begin{align}\label{Jhyouka0924}
J_s(a,b)>0,\ {\partial}_aJ_s(a,b)>0,\ {\partial}_bJ_s(a,b)>0.
\end{align}
By \eqref{Jhyouka0924}, we obtain \eqref{mathcalNsei}. In addition, by Euler's identity for homogeneous functions, we obtain
\begin{align}\label{Jsnohyouka1}
0\leq b{\partial}_bJ_s\leq a{\partial}_aJ_s+b{\partial}_bJ_s=(p+1)J_s.
\end{align}
By \eqref{qes} and \eqref{Jsnohyouka1}, we obtain 
\begin{align*}
\left|\mathcal{N}_s^{\prime}(r)\right|&\lesssim \int {\partial}_bJ_s(Q,Q(\cdot-re_1))Q(\cdot-re_1)dx
\\
&\lesssim \int J_s(Q,Q(\cdot-re_1))dx\lesssim \mathcal{N}_s,
\end{align*}
which implies \eqref{mathcalNbibunes}. Finally, we prove \eqref{mathcalNgensui}. For $0\leq a,b\leq \|Q\|_{L^{\infty}}$, Taylor expansion yields
\begin{align*}
J_s(a,b)\lesssim \left(\max{(a,b)}\right)^{p-1}\left(\min{(a,b)}\right)^2\lesssim (ab)^{\min{(2,\frac{p+1}{2})}}.
\end{align*}
Therefore, we obtain
\begin{align*}
\mathcal{N}_s(r)\lesssim \int \left(QQ(\cdot-re_1)\right)^{\min{(2,\frac{p+1}{2})}}\lesssim e^{-\gamma r},
\end{align*}
which implies \eqref{mathcalNgensui}.
\end{proof}

Furthermore, to estimate interactions involving three or more solitons, we define $\mathcal{R}$ for $z=(z_0,z_1,\ldots,z_k)\in(\mathbb{R}^2)^{k+1}$ by
\begin{align}\label{mathcalRdef}
\mathcal{R}(z)=\frac{1}{p+1}\int\left(|Q_{\sum}|^{p+1}+(k-1)\sum_{i=0}^k|Q_i|^{p+1}-\sum_{0\leq i<j\leq k}|Q_i+Q_j|^{p+1}\right)dx.
\end{align}
Then, we have the following estimate.
\begin{lemma}\label{santailem}
For every $\sqrt{3}<\beta<\sqrt{3}+(2-\sqrt{3})\min{(p-2,1)}$, we have for $D\gg 1$
\begin{align}\label{mathcalRes}
\left|\mathcal{R}(z)\right|+\left|\nabla \mathcal{R}(z)\right|\lesssim e^{-\beta D}.
\end{align}
\end{lemma}
\begin{proof}
First, for $a=(a_0,a_1,\ldots, a_k)\in\mathbb{R}^{k+1}$, we define
\begin{align*}
T(a)=\left|\sum_{i=0}^ka_i\right|^{p+1}+(k-1)\sum_{i=0}^k|a_i|^{p+1}-\sum_{0\leq i<j\leq k}\left|a_i+a_j\right|^{p+1}.
\end{align*}
Furthermore, for $i=0,1,\ldots, k$, we define 
\begin{align*}
h_i(x)=(p+1)\left\{|a_i+x|^{p-1}(a_i+x)-|a_i|^{p-1}a_i\right\}.
\end{align*}
Set $A=\max_{0\leq j\leq k}|a_j|$. Then $h_i(0)=0$ and
\begin{align*}
\sup_{|x|\leq kA}|h_i''(x)|\lesssim A^{p-2}.
\end{align*}
Consequently, we have
\begin{align*}
\left|{\partial}_iT\right|=\left|h_i\left(\sum_{j\neq i}a_j\right)-\sum_{j\neq i}h_i(a_j)\right|\lesssim \left(\max_{0\leq j\leq k}|a_j|\right)^{p-2}\left(\sum_{0\leq j<l\leq k,\ j,l\neq i}|a_ja_l|\right).
\end{align*}
By Euler's identity for homogeneous functions, we obtain
\begin{align}\label{Thyoukaganbatta}
|T|+\sum_{i=0}^k\left|a_i{\partial}_iT\right|\lesssim \max_{i\neq j,j\neq l, l\neq i} |a_i|^{p-1}|a_j||a_l|.
\end{align}
By \eqref{qes} and \eqref{Thyoukaganbatta}, we have
\begin{align}\label{qhyouka}
|\mathcal{R}(z)|+|\nabla \mathcal{R}(z)|\lesssim \int \max_{i\neq j,j\neq l, l\neq i} |Q_i|^{p-1}|Q_j||Q_l|.
\end{align}
Here, we claim that for every $x\in\mathbb{R}^2$ and $i,j,l\in\{0,1,\ldots, k\}$ satisfying $i\neq j, j\neq l, l\neq i$, we have
\begin{align}\label{nijigenkikaes}
(p-1)|x-z_i|+|x-z_j|+|x-z_l|\geq \left\{ \sqrt{3}+(2-\sqrt{3})\min{(p-2,1)} \right\} D.
\end{align}
Indeed, if the three points $z_i,z_j,z_l$ in the plane are collinear or form a triangle with an angle greater than $120^\circ$, it is well known that the sum of the distances to these points attains its minimum at one of the three points. Consequently, for every $x\in\mathbb{R}^2$, we have
\begin{align*}
|x-z_i|+|x-z_j|+|x-z_l|\geq 2D.
\end{align*}
Furthermore, if the three points form a triangle whose interior angles are all at most $120^\circ$, it is well known that the minimizing point $x$ is the Fermat point. In this case, setting $m_1=\min{(|x-z_i|,|x-z_j|,|x-z_l|)}$ and $m_2=\m{(|x-z_i|,|x-z_j|,|x-z_l|)}$, we obtain
\begin{align*}
D^2\leq m_1^2+m_1m_2+m_2^2\leq \frac{\left(|x-z_i|+|x-z_j|+|x-z_l|\right)^2}{3}.
\end{align*}
Thus, in either case, we obtain
\begin{align}\label{fermattenes}
|x-z_i|+|x-z_j|+|x-z_l|\geq \sqrt{3}D.
\end{align}
For further details, see \cite{MN2}. Furthermore, we have
\begin{align}\label{sankakuhutou}
2D&\leq |z_i-z_j|+|z_i-z_l|\leq 2|x-z_i|+|x-z_j|+|x-z_l|.
\end{align}
By \eqref{fermattenes} and \eqref{sankakuhutou}, we obtain \eqref{nijigenkikaes}. Here, we define
\begin{align*}
s_p(x)=(p-1)|x-z_i|+|x-z_j|+|x-z_l|.
\end{align*}
Then, by $|Q(x)|\lesssim e^{-|x|}$, we have $|Q_i(x)|^{p-1}|Q_j(x)||Q_l(x)|\lesssim e^{-s_p(x)}$. Therefore, by \eqref{nijigenkikaes}, for every $\sqrt{3}<\beta<\sqrt{3}+(2-\sqrt{3})\min{(p-2,1)}$, we obtain
\begin{align}\label{fermatinteraction}
\begin{aligned}
\int_{\mathbb{R}^2}|Q_i|^{p-1}|Q_j||Q_l|&\lesssim \int_{\mathbb{R}^2}e^{-s_p(x)}dx
\\
&\lesssim \int_{|x-z_i|\leq 2D}e^{-s_p(x)}dx+\int_{|x-z_i|>2D}e^{-|x-z_i|}dx
\\
&\lesssim D^2e^{- \left\{ \sqrt{3}+(2-\sqrt{3})\min{(p-2,1)} \right\} D}+De^{-2D}\lesssim e^{-\beta D}.
\end{aligned}
\end{align}
By \eqref{qhyouka} and \eqref{fermatinteraction}, we obtain \eqref{mathcalRes}.

\end{proof}

For the subsequent analysis, we fix a constant $\beta_{\star}$ satisfying
\begin{align}\label{betastardef}
\sqrt{3}<\beta_{\star}<\min{(\frac{87}{50}, \sqrt{3}+(2-\sqrt{3})\min{(p-2,1)})}<\frac{7}{4}.
\end{align}

\section{Refined dynamics of the soliton centers}

In this section, we refine the analysis of the dynamics of the soliton centers using improved interaction estimates that exploit the two-dimensional geometry. First, we define
\begin{align}\label{Zdef}
\begin{aligned}
&Z_i=z_i-z_0,\ \rho_i=\left|Z_i\right|\ \ (1\leq i\leq k),
\\
&Z_{ij}=z_j-z_i,\ \rho_{ij}=\left|Z_{ij}\right|\ \ (1\leq i,j\leq k,\ i\neq j).
\end{aligned}
\end{align}
Furthermore, we define the following quantity associated with the energy perturbation:
\begin{align}\label{Gdef}
G(z)=\sum_{i=1}^k\Psi(\rho_i)-\sum_{1\leq i<j\leq k}\Psi(\rho_{ij})-\sum_{i=1}^k\mathcal{N}_{-1}(\rho_i)-\sum_{1\leq i<j\leq k}\mathcal{N}_{1}(\rho_{ij}).
\end{align}
The function $G$ represents the leading-order contribution to the energy difference $E(\vec{u})-(k+1)E(Q,0)$. Indeed, for $(1,k)$-sign soliton solutions, by direct computation and Lemma \ref{vareprefinelem}, we have
\begin{align}\label{Eexpansion}
\begin{aligned}
E(\vec{u})-(k+1)E(Q,0)&=E(\vec{Q}_{\sum}+\vec{\varepsilon})-E(\vec{Q}_{\sum})+E(\vec{Q}_{\sum})-(k+1)E(Q,0)
\\
&=E(\vec{Q}_{\sum})-(k+1)E(Q,0)+O(g(D)^2).
\end{aligned}
\end{align}
Furthermore, when we define
\begin{align}\label{Udef}
\mathcal{U}(z)=E(\vec{Q}_{\sum})-(k+1)E(Q,0),
\end{align}
we have the following lemma.
\begin{lemma}\label{mathcalUmeijilem}
We have for $D\gg 1$ and $0\leq i\leq k$,
\begin{align}\label{U-Ges}
\left|\mathcal{U}-G\right|\lesssim e^{-\beta_{\star}D}, \ \left|\nabla_{z_i} G-\langle R_2,\nabla Q_i\rangle \right|\lesssim e^{-\beta_{\star}D}. 
\end{align}
\end{lemma}
\begin{proof}
By direct computation, we have for $1\leq i\leq k$ and $1\leq i<j\leq k$, respectively,
\begin{align*}
E(Q_0+Q_i,0)-2E(Q,0)&=\Psi(\rho_i)-\mathcal{N}_{-1}(\rho_i),
\\
E(Q_i+Q_j,0)-2E(Q,0)&=-\Psi(\rho_{ij})-\mathcal{N}_1(\rho_{ij}).
\end{align*}
Furthermore, by direct computation and  \eqref{mathcalRdef}, we have
\begin{align}\label{Uequation}
\mathcal{U}(z)=\sum_{0\leq i<j\leq k}\left\{E(Q_i+Q_j,0)-2E(Q,0)\right\}-\mathcal{R}(z)&=G(z)-\mathcal{R}(z).
\end{align}
In addition, since we have \eqref{Udef} and
\begin{align*}
-\Delta Q_{\sum}+Q_{\sum}-f(Q_{\sum})&=\sum_{i=0}^k\left(-\Delta Q_i+Q_i\right)-f(Q_{\sum})=-R_2
\end{align*}
and $\nabla_{z_i}Q_{\sum}=-\nabla_xQ_i$, we obtain
\begin{align}\label{nablaUequation}
\begin{aligned}
\nabla_{z_i}\mathcal{U}&=\langle\nabla_xQ_{\sum},\nabla_x(\nabla_{z_i}Q_{\sum})\rangle+\langle Q_{\sum}-f(Q_{\sum}),\nabla_{z_i}Q_{\sum}\rangle
\\
&=\langle-\Delta Q_{\sum}+Q_{\sum}-f(Q_{\sum}),-\nabla_xQ_i\rangle
\\
&=\langle R_2,\nabla_xQ_i\rangle.
\end{aligned}
\end{align}
By \eqref{mathcalRes}, \eqref{Uequation}, and \eqref{nablaUequation} , we obtain \eqref{U-Ges}.
\end{proof}
We use this estimate to refine our analysis of the dynamics of the soliton centers. Before stating the next lemma, we define
\begin{align}
C_{\star}=2\alpha\|{\partial}_1Q\|_{L^2}^2. \label{Cstardef}
\end{align}

\begin{lemma}\label{chushinseimitulem}
Let $\vec{u}$ be a $(1,k)$-sign soliton solution satisfying \eqref{1,ksolicondrevise}. Then, we have for $t\gg 1$ and for $0\leq i\leq k$,
\begin{align}\label{zesref0927}
C_{\star}\dot{z}_i=-\nabla_{z_i}G(z)+O\left(e^{-\beta_{\star}D}\right).
\end{align}
\end{lemma}
\begin{proof}
We note that by \eqref{prees1} and Lemma \ref{vareprefinelem}, we have
\begin{align}\label{g(D)-1order}
|\dot z|+\|\vec{\varepsilon}\|_{\mathcal{H}}\lesssim g(D).
\end{align}
Fix $0\leq i\leq k$ and $l\in\{1,2\}$. By \eqref{varepnoode}, \eqref{1,ksolicondrevise}, and \eqref{g(D)-1order}, we have
\begin{align}\label{zinotoushiki}
\begin{aligned}
0&=\int \left({\partial}_t\eta+2\alpha{\partial}_t\varepsilon\right){\partial}_lQ_i+\int \left(\eta+2\alpha\varepsilon\right){\partial}_t{\partial}_lQ_i
\\
&=\langle (\Delta-1+f^{\prime}(Q_{\sum}))\varepsilon,{\partial}_lQ_i\rangle+\langle R_2,{\partial}_lQ_i\rangle+2\alpha\langle R_1,{\partial}_lQ_i\rangle
\\
&\quad+\langle R_3,{\partial}_lQ_i\rangle-\langle \eta+2\alpha\varepsilon,\dot{z}_i\cdot \nabla{\partial}_lQ_i\rangle.
\end{aligned}
\end{align}
Furthermore, by \eqref{R2taylores}, we have
\begin{align}\label{Qsenkeikasayousohyouka}
\begin{aligned}
\langle (\Delta-1+f^{\prime}(Q_{\sum}))\varepsilon,{\partial}_lQ_i\rangle&=-\langle(-\Delta+1-f^{\prime}(Q_i))\varepsilon,{\partial}_lQ_i\rangle
\\
&\quad +\langle (f^{\prime}(Q_{\sum})-f^{\prime}(Q_i))\varepsilon,{\partial}_lQ_i\rangle 
\\
&=-\langle(-\Delta+1-f^{\prime}(Q_i)){\partial}_lQ_i,\varepsilon\rangle+O(g(D)\|\vec{\varepsilon}\|_{\mathcal{H}})
\\
&=O\left(g(D)^2\right).
\end{aligned}
\end{align}
In addition, by direct computation, we have
\begin{align}\label{iroiroes0927}
\begin{aligned}
&\langle R_1,{\partial}_lQ_i\rangle=\sum_{j=0}^k\langle \dot{z}_j\cdot \nabla Q_j,{\partial}_lQ_i\rangle=\dot{z}_{i,l}\|{\partial}_lQ_i\|_{L^2}^2+O(g(D)^{\frac{7}{4}}),
\\
&\langle R_3,{\partial}_lQ_i\rangle=O\left(\|\vec{\varepsilon}\|_{\mathcal{H}}^2\right)=O(g(D)^2),
\\
&\langle \eta+2\alpha\varepsilon,\dot{z}_i\cdot \nabla{\partial}_lQ_i\rangle=O(|\dot{z}|\|\vec{\varepsilon}\|_{\mathcal{H}})=O(g(D)^2).
\end{aligned}
\end{align}
We note that the estimate for $R_1$ follows from the bound $\int |\nabla Q_i||\nabla Q_j|\,dx\lesssim e^{-\frac{9}{10}D}$ for $i\neq j$. By \eqref{zinotoushiki}, \eqref{Qsenkeikasayousohyouka}, and \eqref{iroiroes0927}, we obtain
\begin{align}\label{zchottohyoukayosi}
2\alpha\|{\partial}_1Q\|_{L^2}^2\dot{z}_{i,l}=-\langle R_2,{\partial}_lQ_i\rangle+O(g(D)^{\frac{7}{4}}).
\end{align}
By Lemma \ref{mathcalUmeijilem} and \eqref{zchottohyoukayosi}, we obtain \eqref{zesref0927}.
\end{proof}
Furthermore, the relation between $G$ and the energy difference yields the following lemma.

\begin{lemma}\label{Genergy_lower_prop}
Let $\vec{u}$ be a $(1,k)$-sign soliton solution satisfying \eqref{1,ksolicondrevise}. Then, there exists $C_{\circ}>0$ such that, for all sufficiently large $t$,
\begin{align}\label{Genergylower}
G(z(t))\geq-C_{\circ}e^{-\beta_\star D(t)}.
\end{align}
\end{lemma}
\begin{proof}
Since the energy $E(\vec{u}(t))$ is nonincreasing and $\lim_{t\to\infty}E(\vec{u}(t))=(k+1)E(Q,0)$, we have
\begin{align*}
0\leq E(\vec{u}(t))-(k+1)E(Q,0)=\mathcal{U}(t)+O(g(D(t))^2).
\end{align*}
Therefore, we obtain
\begin{align}\label{mathcalUupper}
-\mathcal{U}(t)\lesssim g(D(t))^2.
\end{align}
By \eqref{U-Ges} and \eqref{mathcalUupper}, there exists $C_{\circ}>0$ such that \eqref{Genergylower} holds. Therefore, we complete the proof.
\end{proof}

\section{Hexagonal rigidity and a gradient lower bound}

In this section, we show that, under the lower bound \eqref{Genergylower}, the vanishing of the normalized gradient
forces a regular hexagonal cluster to emerge in the limiting configuration. We first derive the limiting balance equations and then classify the configurations satisfying them. First, we define 
\begin{align}\label{Adef}
\mathcal{A}(z)=\sum_{i=1}^k \Psi(\rho_i).
\end{align}

The aim of this section is to prove the following proposition. For a configuration $z=(z_0,\ldots,z_k)$, we write
\begin{align*}
D=D(z):=\min_{0\leq i<j\leq k}|z_j-z_i|.
\end{align*}

\begin{proposition}\label{haitipropo}
For every $c>0$, there exist $D_1>0$ and $c_1>0$ such that, if $z\in(\mathbb{R}^2)^{k+1}$ satisfies
\begin{align}\label{haitikatei}
D\geq D_1,\ G(z)\geq -ce^{-\beta_{\star}D},
\end{align}
we have
\begin{align}\label{nablaGlower}
\left|\nabla_zG(z)\right|^2\geq c_1\mathcal{A}(z)\left(\max{(G(z),0)}+\Psi(\sqrt{3}D)\right).
\end{align}
\end{proposition}

\subsection{Limiting balance equations and hexagonal rigidity}

In this subsection, we show that, under the lower bound \eqref{Genergylower}, the vanishing of the normalized gradient implies that the limiting configuration contains six points at the vertices of a regular hexagon. We first derive the limiting balance equations and then classify the configurations satisfying them.

\begin{lemma}\label{limitingbalance}
Let $z^{(n)}\in(\mathbb R^2)^{k+1}$ be a sequence satisfying
\begin{align}\label{znassumption}
\begin{aligned}
&\lim_{n\to\infty}D(z^{(n)})=\infty,\ \ \lim_{n\to\infty}\frac{|\nabla_z G(z^{(n)})|}{\mathcal{A}(z^{(n)})}=0,
\\
&G(z^{(n)})\geq -C_{\circ}e^{-\beta_{\star}D}.
\end{aligned}
\end{align}
Then, there exist $(\tilde{z}^{(n)})\subset (z^{(n)})$, $a_i\geq 0,\ \mathfrak{r}_i\in [1,\infty],\ \omega_i\in\mathbb{R}^2$ ($1\leq i\leq k$),\ $b_{ij}\geq 0,\ \mathfrak{r}_{ij}\in[1,\infty],\ \omega_{ij}\in\mathbb{R}^2$ ($1\leq i\neq j\leq k$) such that 
\begin{align}\label{abomegabubun}
\begin{aligned}
&a_i=\lim_{n\to\infty}\frac{\Psi(\rho_i(\tilde{z}^{(n)}))}{\mathcal{A}(\tilde{z}^{(n)})},\ b_{ij}=\lim_{n\to\infty}\frac{\Psi(\rho_{ij}(\tilde{z}^{(n)}))}{\mathcal{A}(\tilde{z}^{(n)})},
\\
&\mathfrak{r}_i=\lim_{n\to\infty}\frac{\rho_i(\tilde{z}^{(n)})}{D(\tilde{z}^{(n)})},\ \mathfrak{r}_{ij}=\lim_{n\to\infty}\frac{\rho_{ij}(\tilde{z}^{(n)})}{D(\tilde{z}^{(n)})},
\\
&\omega_i=\lim_{n\to\infty}\frac{\tilde{z}_i^{(n)}-\tilde{z}_0^{(n)}}{\rho_i(\tilde{z}^{(n)})},\ \omega_{ij}=\lim_{n\to\infty}\frac{\tilde{z}_j^{(n)}-\tilde{z}_i^{(n)}}{\rho_{ij}(\tilde{z}^{(n)})}.
\end{aligned}
\end{align}
Furthermore, for $1\leq i,j\leq k$ with $i\neq j$, we have
\begin{align}\label{omeganonagasa}
|\omega_i|=|\omega_{ij}|=1,\ \omega_{ji}=-\omega_{ij}.
\end{align}
In addition, when we define
\begin{align}\label{Lambdateigi1}
\Lambda=\{i\in\{1,2,\ldots,k\}:\mathfrak{r}_i<\infty\},
\end{align}
we have 
\begin{align}\label{turiaieqiroiro}
\begin{aligned}
&\sum_{i\in \Lambda}a_i=1,\ \ \sum_{i<j,\ i,j\in\Lambda}b_{ij}\leq 1,\ \ \sum_{i\in \Lambda} a_i\omega_i=0,
\\
&a_i\omega_i+\sum_{j\in \Lambda, j\neq i}b_{ij}\omega_{ij}=0\ (i\in\Lambda).
\end{aligned}
\end{align}
In particular, we have
\begin{align}\label{nagasadekaito0}
\begin{aligned}
&a_i=0\ \ (\mathfrak r_i>1),\ b_{ij}=0\quad(\mathfrak r_{ij}>1),\ b_{ij}=0\ \ (i\notin\Lambda\ \text{or}\ j\notin\Lambda).
\end{aligned}
\end{align}
\end{lemma}

\begin{proof}
We suppress the dependence on $n$ when no confusion arises. We define
\begin{align*}
\mathcal{B}
=\sum_{1\leq i<j\leq k}\Psi(\rho_{ij}).
\end{align*}
By \eqref{mathcalNsei}, \eqref{Gdef}, and \eqref{znassumption}, we have
\begin{align}\label{A-Bes}
\mathcal{B}\leq\mathcal{A}+C_{\circ}e^{-\beta_{\star}D}=\mathcal{A}+o(g(D)).
\end{align}
Furthermore, by the definition of $\mathcal{A}$ and $\mathcal{B}$, we obtain
\begin{align}\label{mathcalAsimpsiD}
\mathcal{A}+\mathcal{B}=\sum_{0\leq i<j\leq k}\Psi(|z_j-z_i|)\sim\Psi(D).
\end{align}
By \eqref{A-Bes} and \eqref{mathcalAsimpsiD}, we have
\begin{align}\label{Asimg(D)}
\mathcal{A}\sim g(D).
\end{align}
In addition, by \eqref{Asimg(D)}, for all sufficiently large $n$,
\begin{align*}
0\leq \frac{\Psi(\rho_i)}{\mathcal A},
\frac{\Psi(\rho_{ij})}{\mathcal A}\lesssim1.
\end{align*}
Moreover, we have $\frac{\rho_i}{D},\frac{\rho_{ij}}{D}\geq 1$ and $\frac{Z_i}{\rho_i},\frac{Z_{ij}}{\rho_{ij}}\in S^1$.
Therefore, by passing to a subsequence $\{\tilde z^{(n)}\}$, we obtain all the limits in \eqref{abomegabubun}.
Here, for the distance ratios, we use the compactness of their reciprocals in $[0,1]$, allowing the limits
$\mathfrak r_i$ and $\mathfrak r_{ij}$ to equal $\infty$. Then, by definition, we have $b_{ij}=b_{ji}$, and $|\omega_i|=|\omega_{ij}|=1$, and $\omega_{ji}=-\omega_{ij}$. This proves \eqref{omeganonagasa}. Furthermore, since we have
\begin{align*}
&\sum_{i=1}^k\frac{\Psi(\rho_i(z^{(n)}))}{\mathcal{A}(z^{(n)})}=1,\ \sum_{1\leq i<j\leq k}\frac{\Psi(\rho_{ij}(z^{(n)}))}{\mathcal{A}(z^{(n)})}=\frac{\mathcal{B}(z^{(n)})}{\mathcal{A}(z^{(n)})},
\end{align*}
by \eqref{A-Bes} and \eqref{Asimg(D)}, we obtain
\begin{align}\label{atobes}
\sum_{i=1}^k a_i=1,\ \sum_{1\leq i<j\leq k}b_{ij}\leq1.
\end{align}

Next, by Lemma \ref{mathcalNeslem} and \eqref{Psies0922} and \eqref{Asimg(D)}, we have
\begin{align}\label{znablaGes}
\begin{aligned}
\frac{\nabla_{z_0}G}{\mathcal A}
&=
\sum_{i=1}^k
\frac{g(\rho_i)}{\mathcal A}
\frac{z_i-z_0}{\rho_i}
+o(1),
\\
\frac{\nabla_{z_i}G}{\mathcal A}
&=
-\frac{g(\rho_i)}{\mathcal A}
\frac{z_i-z_0}{\rho_i}
-\sum_{\substack{1\leq j\leq k\\j\neq i}}
\frac{g(\rho_{ij})}{\mathcal A}
\frac{z_j-z_i}{\rho_{ij}}
+o(1).
\end{aligned}
\end{align}
By \eqref{Psies0922}, \eqref{znassumption}, and \eqref{znablaGes}, we obtain
\begin{align}\label{haitieqconv}
\begin{aligned}
&\sum_{i=1}^k a_i\omega_i=0,
\\
&a_i\omega_i+
\sum_{\substack{1\leq j\leq k\\j\neq i}}
b_{ij}\omega_{ij}=0
\qquad(1\leq i\leq k).
\end{aligned}
\end{align}

We next prove \eqref{nagasadekaito0}.
Suppose that $\mathfrak r_i>1$.
Then there exists $\delta>0$ such that
$\rho_i\geq(1+\delta)D$ for all sufficiently large $n$.
By \eqref{Asimg(D)} and the asymptotics of $\Psi$, we have
\begin{align*}
0\leq\frac{\Psi(\rho_i)}{\mathcal A}\lesssim\frac{g((1+\delta)D)}{g(D)}\lesssim e^{-\delta D/2},
\end{align*}
and therefore we have $a_i=0$. The same argument gives $b_{ij}=0$ if $\mathfrak r_{ij}>1$. In particular, $a_i=0$ for $i\notin\Lambda$. We note that by \eqref{atobes}, we have $\sum_{i\in\Lambda}a_i=1$, and hence $\Lambda\neq\varnothing$. In addition, for $i\in\Lambda$ and $j\notin\Lambda$, we have
\begin{align*}
\frac{\rho_{ij}}D\geq\frac{\rho_j}D-\frac{\rho_i}D\to\infty.
\end{align*}
Thus, we have $\mathfrak r_{ij}=\infty$ and $b_{ij}=0$. Here, for $i\in\Lambda$, we define
\begin{align*}
y_i=\mathfrak r_i\omega_i=\lim_{n\to\infty}\frac{\tilde z_i^{(n)}-\tilde z_0^{(n)}}{D(\tilde z^{(n)})}.
\end{align*}
By \eqref{abomegabubun}, for distinct $i,j\in\Lambda$,
we have
\begin{align*}
&|y_i|=\mathfrak r_i,\ |y_j-y_i|=\mathfrak r_{ij},\ y_j-y_i=\mathfrak r_{ij}\omega_{ij},
\\
&a_i\omega_i=a_i y_i,\ b_{ij}\omega_{ij}=b_{ij}(y_j-y_i).
\end{align*}
Since $b_{ij}=0$ for $i\in\Lambda$ and $j\notin\Lambda$, \eqref{haitieqconv} yields, for $i\in\Lambda$,
\begin{align*}
a_i y_i+\sum_{j\in\Lambda, j\neq i}b_{ij}(y_j-y_i)=0.
\end{align*}
Taking the scalar product with $y_i$ and summing over $i\in\Lambda$, we obtain
\begin{align*}
0=\sum_{i\in\Lambda}a_i|y_i|^2-\sum_{i<j,\ i,j\in\Lambda}b_{ij}|y_j-y_i|^2=1-\sum_{i<j,\ i,j\in\Lambda}b_{ij}.
\end{align*}
Consequently, by \eqref{atobes},
\begin{align*}
1=\sum_{i<j,\ i,j\in\Lambda}b_{ij}\leq \sum_{1\leq i<j\leq k}b_{ij}\leq1.
\end{align*}
Since all $b_{ij}$ are nonnegative, every coefficient $b_{ij}$ with at least one index outside $\Lambda$
must vanish. This proves \eqref{nagasadekaito0}. Finally, restricting \eqref{atobes} and \eqref{haitieqconv}
to $\Lambda$ yields \eqref{turiaieqiroiro}.
\end{proof}
Next, we show that any configuration satisfying \eqref{turiaieqiroiro} and \eqref{nagasadekaito0} contains six points forming a regular hexagon centered at the origin, while the remaining $k-6$ points do not contribute to the limiting
balance equations. To state the next lemma, for $i=1,2,\ldots,6$, we define
\begin{align*}
u_i=\left(\cos{\frac{(i-1)\pi}{3}},\sin{\frac{(i-1)\pi}{3}}\right).
\end{align*}

\begin{lemma}\label{hexagonclassification}
Let $1\leq k^{\prime}\leq k$. We assume that $y_i\in\mathbb{R}^2,\ a_i\geq 0$ ($1\leq i\leq k^{\prime}$) and $b_{ij}\geq 0$ ($i\neq j,\ 1\leq i,j\leq k^{\prime}$) satisfy
\begin{align}
&|y_i|\geq1,\ |y_i-y_j|\geq1\ \ (i\neq j), \label{normalizedseparation}
\\
&b_{ij}=b_{ji}, \label{bij=bji}
\\
&\sum_{i=1}^{k^{\prime}}a_i=1,\ a_i=0\ \text{if}\ |y_i|>1,\ b_{ij}=0\ \text{if}\ |y_i-y_j|>1, \label{contactcoefficients}
\\
&a_i y_i+\sum_{j\neq i}b_{ij}(y_j-y_i)=0\ \ (1\leq i\leq k^{\prime}). \label{zerobalance}
\end{align}
Then, we have $k^{\prime}\geq6$. Furthermore, there exist $\{i_1,i_2,\ldots, i_{k^{\prime}}\}=\{1,2,\ldots,k^{\prime}\}$ and $R\in SO(2)$ such that 
\begin{align}\label{hexagonalpositions}
y_{i_j}=R\left(\cos\frac{(j-1)\pi}{3},\sin\frac{(j-1)\pi}{3}\right)\ \ (1\leq j\leq6),
\end{align}
and 
\begin{align}\label{hexagonalweights}
\begin{aligned}
a_{i_j}&=
\begin{cases}
\frac{1}{6}\ \ &(1\leq j\leq6),
\\
0\ \ &(7\leq j\leq k^{\prime}),
\end{cases}
\\
b_{i_ji_l}&=
\begin{cases}
\frac{1}{6}\ \ &(j,l\in\{1,\ldots,6\},\ j-l\equiv\pm1\pmod6),
\\
0\ \ &(\text{otherwise}),
\end{cases}
\end{aligned}
\end{align}
and
\begin{align}\label{remainingradius}
|y_{i_j}|\geq\sqrt3\ \ (7\leq j\leq k^{\prime}).
\end{align}
\end{lemma}
\begin{proof}
We divide the proof into three steps.

\textbf{Step 1. The points with nonzero coefficients lie on
the unit circle.}

We define
\begin{align*}
I=\left\{ i\in\{1,\ldots,k^{\prime}\}:a_i+\sum_{1\leq j\leq k^{\prime},\ j\neq i}b_{ij}>0\right\}.
\end{align*}
By \eqref{contactcoefficients}, we have $I\neq\varnothing$. We choose $i\in I$ such that $|y_i|=\max_{j\in I}|y_j|$. Furthermore, we assume that $|y_i|>1$. Then, by \eqref{contactcoefficients}, we have $a_i=0$, and hence there exists $j$ such that $b_{ij}>0$. For every such $j$, \eqref{bij=bji} implies $j\in I$, and \eqref{normalizedseparation} and \eqref{contactcoefficients} give $|y_j-y_i|=1$. Therefore, we obtain
\begin{align}\label{es1011}
(y_j-y_i)\cdot y_i=\frac{|y_j|^2-|y_i|^2-|y_j-y_i|^2}{2}\leq-\frac{1}{2}.
\end{align}
By \eqref{zerobalance} and \eqref{es1011}, we obtain
\begin{align*}
0=\sum_{j\neq i}b_{ij}(y_j-y_i)\cdot y_i\leq-\frac{1}{2}\sum_{j\neq i}b_{ij}<0,
\end{align*}
which is a contradiction. By this argument and \eqref{normalizedseparation}, we obtain
\begin{align*}
|y_i|=1\ \ (i\in I).
\end{align*}
In addition, if $b_{ij}>0$, then $i,j\in I$, and therefore we obtain
\begin{align}\label{seisankaku1005}
|y_i|=|y_j|=|y_j-y_i|=1,\ y_i\cdot y_j=\frac{1}{2}.
\end{align}
Thus, taking the scalar product of \eqref{zerobalance} with $y_i$ and using \eqref{seisankaku1005}, we obtain
\begin{align*}
a_i=\frac{1}{2}\sum_{j\neq i}b_{ij}\ \ (i\in I).
\end{align*}
In particular, the definition of $I$ gives
\begin{align*}
I=\{i\in\{1,\ldots,k^{\prime}\}:a_i>0\}.
\end{align*}

\textbf{Step 2. The points indexed by $I$ form a regular hexagon.}

All points in $\{y_i:i\in I\}$ lie on the unit circle, and the distance between any two distinct points is at least $1$.
Thus, their successive angular gaps are at least $\frac{\pi}{3}$, and consequently we have $|I|\leq6$.

Fix $i\in I$. By Step 1, at least one $b_{ij}$ is positive, and $b_{ij}>0$ implies $|y_j|=|y_j-y_i|=1$. Thus, at most two indices $j\neq i$ can have $b_{ij}>0$, corresponding to the two possible neighboring points at angles $\pm\frac{\pi}{3}$ from $y_i$. Then, \eqref{zerobalance} implies that both neighboring points, denoted by $y_j$ and $y_l$, are present and that $a_i=b_{ij}=b_{il}$. Repeating this argument at the neighboring points produces six equally spaced points on the unit circle. Since $|I|\leq6$, these are precisely the points indexed by $I$. In particular, we have $k^{\prime}\geq6$.

Here, we choose their indices $i_1,\ldots,i_6$ in counterclockwise order, and enumerate the remaining indices arbitrarily as
$i_7,\ldots,i_{k^{\prime}}$. Let $R\in SO(2)$ be the rotation satisfying $R(1,0)=y_{i_1}$. Since the six points are equally spaced, we have
\begin{align*}
y_{i_j}=R\left(\cos\frac{(j-1)\pi}{3},\sin\frac{(j-1)\pi}{3}\right)\ \ (1\leq j\leq6).
\end{align*}
This proves \eqref{hexagonalpositions}. Then, noting that, by the argument in Step 1, all coefficients $b_{ij}$ other than those corresponding to the two neighboring points of $y_i$ vanish, we obtain \eqref{hexagonalweights}
by considering \eqref{zerobalance} for $i=i_1,i_2,\ldots,i_6$ and using the normalization $\sum_{i=1}^{k^{\prime}}a_i=1$ in \eqref{contactcoefficients}.

\textbf{Step 3. Location of the remaining points.}

Fix $7\leq j\leq k^{\prime}$ and set $r=|y_{i_j}|\geq1$. Then, there exists $1\leq l\leq6$ such that the angle between $y_{i_j}$ and $y_{i_l}$ is at most $\frac{\pi}{6}$. By \eqref{normalizedseparation}, we have
\begin{align*}
1\leq |y_{i_j}-y_{i_l}|^2\leq r^2+1-2r\cos\frac{\pi}{6}=r^2+1-\sqrt{3}r.
\end{align*}
Since we have $r\geq1$, we obtain $r\geq\sqrt{3}$, which proves \eqref{remainingradius}.
\end{proof}

\subsection{Estimates for deviations from a regular hexagon}

In this subsection, we derive estimates for the deviation of a configuration from a regular hexagon. First, we define $H_0:(\mathbb{R}^2)^7\to \mathbb{R}$ as 
\begin{align*}
H_0(z_0,z_1,z_2,\ldots,z_6)=\sum_{i=1}^6\Psi(|z_i-z_0|)-\sum_{i=1}^5\Psi(|z_{i+1}-z_i|)-\Psi(|z_6-z_1|).
\end{align*}
We note that if $z_1,z_2,\ldots,z_6$ form, in this order, the vertices of a regular hexagon centered at $z_0$, then $H_0(z_0,\ldots,z_6)=0$. In addition, for $z=(z_0,z_1,\ldots,z_6)$, let $(c,r,R)$ be a minimizer of the following minimization problem:
\begin{align*}
\min_{c\in\mathbb{R}^2,\ r\geq0,\ R\in SO(2)}\left\{|z_0-c|^2+\sum_{i=1}^6|z_i-c-r(Ru_i)|^2\right\}.
\end{align*}
We then define $q=(q_0,q_1,\ldots,q_6)\in(\mathbb{R}^2)^7$ by
\begin{align*}
q_0=c,\ q_i=c+r(Ru_i)\ \ (1\leq i\leq6).
\end{align*}
Furthermore, we define $w=(w_0,w_1,\ldots, w_6)\in(\mathbb{R}^2)^7$ as
\begin{align*}
w=z-q.
\end{align*}
Then, by the minimality of $|w|^2$, we have
\begin{align}\label{wchokkoujouken}
\sum_{i=0}^6w_i=0,\ \sum_{i=1}^6w_i\cdot u_i=0,\ \sum_{i=1}^6w_i\cdot Ju_i=0,
\end{align}
where $J:\mathbb{R}^2\to\mathbb{R}^2$ is defined by
\begin{align*}
J(x,y)=(-y,x).
\end{align*}
We then obtain the following estimates near a regular hexagonal configuration.
\begin{lemma}\label{H0nogosalemma}
There exist $r_H,\delta_H,C_H,c_H>0$ such that if $z=(z_0,\ldots,z_6)\in(\mathbb R^2)^7$ satisfies
\begin{align*}
r=|q_1-q_0|\geq r_H,\ |w|\leq\delta_H,
\end{align*}
then
\begin{align}\label{H0nogosahyouka}
|H_0(z)|\leq C_Hg(r)|w|^2,\ |\nabla_zH_0(z)|\geq c_Hg(r)|w|.
\end{align}
\end{lemma}

\begin{proof}
By translation and rotation invariance, we may assume that
\begin{align*}
q_0=0,\ \ q_i=ru_i\ \ (1\leq i\leq6).
\end{align*}
These transformations preserve the orthogonality conditions \eqref{wchokkoujouken}. Throughout the proof, the outer indices are interpreted cyclically modulo $6$. Since all six radii and all six adjacent edges have length $r$, and we have
\begin{align*}
\sum_{i=1}^6u_i=0,\ u_{i-1}+u_{i+1}=u_i,
\end{align*}
we obtain $H_0(q)=0$ and $\nabla_zH_0(q)=0$. We define $v_i=Ju_i$. By \eqref{wchokkoujouken}, there exist $\mathfrak{s}_i,\mathfrak{t}_i\in\mathbb{R}$ ($1\leq i\leq6$) such that
\begin{align}\label{wdefalphabeta}
\begin{aligned}
w_0&=-\frac{1}{7}\sum_{i=1}^6(\mathfrak{s}_iu_i+\mathfrak{t}_iv_i),
\\
w_i&=w_0+\mathfrak{s}_iu_i+\mathfrak{t}_iv_i\ \ (1\leq i\leq 6),
\\
\sum_{i=1}^6\mathfrak{s}_i&=\sum_{i=1}^6\mathfrak{t}_i=0.
\end{aligned}
\end{align}
Indeed, for each $i=1,\ldots,6$, the vector $w_i-w_0$ can be expressed as a linear combination of the orthonormal basis
vectors $u_i$ and $v_i$. Thus, defining
\begin{align*}
\mathfrak{s}_i=(w_i-w_0)\cdot u_i,\ \ \mathfrak{t}_i=(w_i-w_0)\cdot v_i,
\end{align*}
we have
\begin{align*}
w_i-w_0=\mathfrak{s}_iu_i+\mathfrak{t}_iv_i.
\end{align*}
Furthermore, since we have $\sum_{i=1}^6u_i=\sum_{i=1}^6v_i=0$, we obtain
\begin{align*}
\sum_{i=1}^6\mathfrak{s}_i&=\sum_{i=1}^6w_i\cdot u_i-w_0\cdot \sum_{i=1}^6u_i=0,
\\
\sum_{i=1}^6\mathfrak{t}_i&=\sum_{i=1}^6w_i\cdot v_i-w_0\cdot \sum_{i=1}^6 v_i=0.
\end{align*}
By $\sum_{i=0}^6w_i=0$, we have
\begin{align*}
0=7w_0+\sum_{i=1}^6(\mathfrak{s}_iu_i+\mathfrak{t}_iv_i),
\end{align*}
and we obtain \eqref{wdefalphabeta}.

We now regard $\mathfrak{s}_i,\mathfrak{t}_i$ as variables and set
\begin{align*}
X=\left\{ \xi=(\mathfrak{s},\mathfrak{t})\in\mathbb R^6\times\mathbb R^6: \sum_{i=1}^6\mathfrak{s}_i=\sum_{i=1}^6\mathfrak{t}_i=0
\right\}.
\end{align*}
Furthermore, for $\xi=(\mathfrak{s},\mathfrak{t})\in\mathbb R^{12}$, we define
\begin{align*}
d_i(\xi)&=\mathfrak{s}_iu_i+\mathfrak{t}_iv_i,
\\
W_0(\xi)&=-\frac{1}{7}\sum_{j=1}^6d_j(\xi),
\\
W_i(\xi)&=W_0(\xi)+d_i(\xi)\ \ (1\leq i\leq6).
\end{align*}
To avoid confusion in notation, we denote by $\xi_{\ast}\in X$ the coefficient vector corresponding to the fixed displacement $w$, so that $w=W(\xi_{\ast})$. In what follows, we let $\xi$ vary in $X$, while keeping $\xi_{\ast}$ fixed.

We first compare the norms of $W(\xi)$ and $\xi$. By direct computation, we have
\begin{align*}
|W(\xi)|^2=|\xi|^2-\frac{1}{7}\left|\sum_{i=1}^6d_i(\xi)\right|^2,
\end{align*}
and by the Cauchy--Schwarz inequality, we have 
\begin{align}\label{xisim}
\frac{1}{7}|\xi|^2\leq |W(\xi)|^2\leq|\xi|^2.
\end{align}
For fixed $r$, we define
\begin{align*}
F_r(\xi)=H_0(q+W(\xi)).
\end{align*}
We equip $X$ with the Euclidean inner product and denote the gradient and Hessian of the restriction to $X$ by $\nabla_X$ and $D_X^2$, respectively. Since $H_0(q)=0$ and $\nabla_zH_0(q)=0$, we have $F_r(0)=0$ and $\nabla_XF_r(0)=0$.
 To compute the Hessian, we set $e_i=u_{i+1}-u_i$ and define
\begin{align*}
B_i(\xi)=\frac{\mathfrak{s}_i+\mathfrak{s}_{i+1}}2+\frac{\sqrt3}{2}(\mathfrak{t}_{i+1}-\mathfrak{t}_i),\ C_i(\xi)=\frac{\sqrt3}{2}
(\mathfrak{s}_i-\mathfrak{s}_{i+1})+\frac{\mathfrak{t}_i+\mathfrak{t}_{i+1}}2.
\end{align*}
We note that $B_i(\xi)$ and $C_i(\xi)$ are the components of the edge perturbation $d_{i+1}(\xi)-d_i(\xi)$ along $e_i$ and $Je_i$, respectively. Then, we have
\begin{align*}
ru_i+d_i(\xi)
&=(r+\mathfrak{s}_i)u_i+\mathfrak{t}_iv_i,
\\
re_i+d_{i+1}(\xi)-d_i(\xi)
&=(r+B_i(\xi))e_i+C_i(\xi)Je_i.
\end{align*}
Since both $(u_i,v_i)$ and $(e_i,Je_i)$ are orthonormal
bases, defining
\begin{align*}
f_r(a,b)=\Psi\left(\sqrt{(r+a)^2+b^2}\right),
\end{align*}
we obtain
\begin{align*}
F_r(\xi)
=\sum_{i=1}^6f_r(\mathfrak{s}_i,\mathfrak{t}_i)
-\sum_{i=1}^6f_r(B_i(\xi),C_i(\xi)).
\end{align*}
Furthermore, by direct differentiation, we have
\begin{align*}
D^2f_r(0,0)=
\begin{pmatrix}
\Psi^{\prime \prime}(r)&0\\
0&\frac{\Psi^{\prime}(r)}{r}
\end{pmatrix}.
\end{align*}
Consequently, for every $\xi\in X$, we have
\begin{align*}
D_X^2F_r(0)[\xi,\xi]
&=\left.\frac{d^2}{d\lambda^2}F_r(\lambda\xi)\right|_{\lambda=0}
\\
&=\Psi^{\prime \prime}(r)\sum_{i=1}^6\left(\mathfrak{s}_i^2-B_i(\xi)^2\right)+\frac{\Psi^{\prime}(r)}{r}\sum_{i=1}^6\left(\mathfrak{t}_i^2-C_i(\xi)^2\right).
\end{align*}
Next, we define the quadratic form
\begin{align*}
\mathcal Q(\xi)=\frac{1}{2}\sum_{i=1}^6\left(\mathfrak{s}_i^2-B_i(\xi)^2\right).
\end{align*}
Since $B_i$ and $C_i$ are linear forms with coefficients independent of $r$, \eqref{Psies0922} and the above formula
for $D_X^2F_r(0)$ yield
\begin{align}\label{FrHessianapprox}
\left|D_X^2F_r(0)\xi-g(r)D_X^2\mathcal{Q}\,\xi\right|\lesssim \frac{g(r)}r|\xi|\ \ (\xi\in X).
\end{align}
We next show that $D_X^2\mathcal Q$ has trivial kernel. The partial derivatives of $\mathcal Q$ are
\begin{align*}
\frac{\partial\mathcal Q}{\partial\mathfrak{s}_i}=\mathfrak{s}_i-\frac{B_i+B_{i-1}}2,\ \  \frac{\partial\mathcal Q}{\partial\mathfrak{t}_i}=\frac{\sqrt{3}}{2}(B_i-B_{i-1}),
\end{align*}
where we suppress the dependence of $B_i$ on $\xi$. Furthermore, since we have $\sum_{i=1}^6B_i=\sum_{i=1}^6\mathfrak{s}_i$, we obtain
\begin{align*}
\sum_{i=1}^6\frac{\partial\mathcal Q}{\partial\mathfrak{s}_i}=\sum_{i=1}^6\frac{\partial\mathcal Q}{\partial\mathfrak{t}_i}=0.
\end{align*}
Thus $\nabla\mathcal Q(\xi)\in X$, and hence we have for $\xi\in X$,
\begin{align*}
\nabla_X\mathcal Q(\xi)=\nabla\mathcal Q(\xi).
\end{align*}
Here, we assume that $\xi\in X$ satisfies $\nabla_X\mathcal Q(\xi)=0$. Then, the derivatives with respect to $\mathfrak{t}_i$ imply that all $B_i$ have a common value $b$. Then, the derivatives with respect to $\mathfrak{s}_i$ give $\mathfrak{s}_i=b$ for every $i$. Furthermore, since we have $\sum_i\mathfrak{s}_i=0$, we obtain $b=0$, and therefore we have
\begin{align*}
\mathfrak{s}_i=B_i=0\ \ (1\leq i\leq6).
\end{align*}
In addition, by the definition of $B_i$, it follows that $\mathfrak{t}_{i+1}=\mathfrak{t}_i$ for every $i$. Together with $\sum_i\mathfrak{t}_i=0$, we have for every $1\leq i\leq 6$, $\mathfrak{t}_i=0$. Thus, we have $\xi=0$, which implies that $D_X^2\mathcal Q$ has trivial kernel.

We now use this nondegeneracy and Taylor's formula to obtain a quadratic upper bound for $|F_r(\xi)|$ and a linear lower
bound for $|\nabla_XF_r(\xi)|$ for large $r$ and small $|\xi|$. Since $\mathcal Q$ is a quadratic form, we have $\nabla_X\mathcal Q(\xi)=D_X^2\mathcal Q\,\xi$. By compactness of the unit sphere in $X$, there exists $\kappa>0$ such that
\begin{align*}
|\nabla_X\mathcal Q(\xi)|\geq\kappa|\xi|\ \ (\xi\in X).
\end{align*}
For $|\xi|\leq1$, all distances appearing in $F_r(\xi)$ are $r+O(1)$. Thus, by \eqref{gnozenkin1} and \eqref{Psies0922}, we have for all sufficiently large $r$,
\begin{align}\label{Frderivativebounds}
\sup_{\xi\in X,\ |\xi|\leq1}\left(\|D_X^2F_r(\xi)\|_{\mathrm{op}}+\|D_X^3F_r(\xi)\|_{\mathrm{op}}\right)\lesssim g(r).
\end{align}
By Taylor's formula, \eqref{FrHessianapprox}, and \eqref{Frderivativebounds}, using $F_r(0)=0$ and $\nabla_XF_r(0)=0$, we obtain, for $\xi\in X$ with $|\xi|\leq1$,
\begin{align}\label{FrTaylorestimates}
|F_r(\xi)|\lesssim g(r)|\xi|^2,\ \ \left|\nabla_XF_r(\xi)-g(r)\nabla_X\mathcal Q(\xi)\right|\lesssim g(r)\left(r^{-1}|\xi|+|\xi|^2\right).
\end{align}
Consequently, we obtain
\begin{align*}
|\nabla_XF_r(\xi)|\geq g(r)|\xi|\left(\kappa-\frac Cr-C|\xi|\right).
\end{align*}
We now apply these estimates to $\xi_{\ast}$. By \eqref{xisim}, we have
\begin{align*}
|w|\leq|\xi_*|\leq\sqrt{7}|w|.
\end{align*}
Choosing $r_H$ sufficiently large and $\delta_H>0$ sufficiently small, we have
\begin{align}\label{es1011-2}
|\xi_{\ast}|\leq 1,\ \ |\nabla_XF_r(\xi_{\ast})|\geq\frac{\kappa}{2}g(r)|\xi_{\ast}|
\end{align}
whenever $r\geq r_H$ and $|w|\leq\delta_H$. Finally, we estimate $H_0$ and $\nabla_zH_0$. By \eqref{xisim}, the linear map $W|_X$ has operator norm at most $1$. Since $F_r(\xi)=H_0(q+W(\xi))$, the chain rule yields
\begin{align}\label{es1011-3}
|\nabla_XF_r(\xi_*)|\leq|\nabla_zH_0(q+w)|.
\end{align}
By \eqref{es1011-2} and \eqref{es1011-3}, we have
\begin{align*}
|\nabla_zH_0(z)|\geq\frac{\kappa}{2}g(r)|\xi_{\ast}|\geq\frac{\kappa}{2}g(r)|w|.
\end{align*}
Moreover, by \eqref{xisim} and \eqref{FrTaylorestimates}, we have
\begin{align*}
|H_0(z)|=|F_r(\xi_{\ast})|\lesssim g(r)|\xi_{\ast}|^2\lesssim g(r)|w|^2,
\end{align*}
and we complete the proof.
\end{proof}

Using these lemmas, we prove Proposition \ref{haitipropo}.

\begin{proof}[Proof of Proposition \ref{haitipropo}]
Fix $c>0$ and suppose that the proposition is false. Then, there exists a sequence $z^{(n)}\in(\mathbb R^2)^{k+1}$ such that
\begin{align}\label{prop4.1nokatei}
\begin{aligned}
\lim_{n\to\infty}D(z^{(n)})=\infty,\ G(z^{(n)})\geq-ce^{-\beta_\star D_n},
\\
\lim_{n\to\infty}\frac{|\nabla_zG(z^{(n)})|^2}{\mathcal A(z^{(n)})\left(\max(G(z^{(n)}),0)+\Psi(\sqrt3D_n)\right)}=0.
\end{aligned}
\end{align}
We define $D_n=D(z^{(n)})$ and suppress the dependence on $n$ when no confusion arises. We divide the proof into four steps.

\textbf{Step 1. Extraction and classification of the limiting configuration.}

By \eqref{mathcalNsei} and \eqref{Gdef}, we have
\begin{align*}
\sum_{1\leq i<j\leq k}\Psi(\rho_{ij})\leq\mathcal A+ce^{-\beta_\star D}.
\end{align*}
In addition, by the definition of $D$, we have
\begin{align*}
\Psi(D)\leq\mathcal A+\sum_{1\leq i<j\leq k}\Psi(\rho_{ij})\leq2\mathcal A+ce^{-\beta_\star D}.
\end{align*}
Together with $\mathcal A\leq k\Psi(D)$ and $e^{-\beta_\star D}=o(\Psi(D))$, we obtain
\begin{align}\label{Apsigsim}
\mathcal A\sim\Psi(D)\sim g(D).
\end{align}
Moreover, by $G\leq\mathcal A$ and $\Psi(\sqrt3D)\lesssim\mathcal A$, we obtain
\begin{align}\label{nablaAnokatei}
\lim_{n\to\infty}\frac{|\nabla_zG|}{\mathcal A}=0.
\end{align}
By \eqref{prop4.1nokatei}, \eqref{nablaAnokatei}, and Lemma \ref{limitingbalance}, after passing to a subsequence
and relabeling $z_1,\ldots,z_k$, we may assume that there exist $a_i\geq0$, $\mathfrak{r}_i\in[1,\infty]$, $\omega_i\in\mathbb R^2$ ($1\leq i\leq k$), and $b_{ij}\geq0$, $\mathfrak{r}_{ij}\in[1,\infty]$,
$\omega_{ij}\in\mathbb R^2$ ($1\leq i\neq j\leq k$) satisfying \eqref{abomegabubun}--\eqref{nagasadekaito0}. Furthermore, we may assume that $\Lambda=\{1,2,\ldots,k'\}$ for some $1\leq k'\leq k$.

Here, for $1\leq i\leq k^{\prime}$, we define
\begin{align*}
y_i=\lim_{n\to\infty}\frac{z_i^{(n)}-z_0^{(n)}}{D_n}.
\end{align*}
Then the points $y_i$ and the coefficients $a_i,b_{ij}$ satisfy \eqref{normalizedseparation}--\eqref{zerobalance}.
Thus, by Lemma \ref{hexagonclassification}, we have $k^{\prime}\geq6$. After relabeling $z_1,\ldots,z_{k'}$ once more,
we may assume that there exists $R_{\infty}\in SO(2)$ such that
\begin{align*}
y_i=R_{\infty} u_i\ \ (1\leq i\leq6),
\end{align*}
and
\begin{align}\label{akbkkeisuu}
\begin{aligned}
a_i&=
\begin{cases}
\frac{1}{6}\ \ &(1\leq i\leq6),
\\
0\ \ &(7\leq i\leq k),
\end{cases}
\\
b_{ij}&=
\begin{cases}
\frac{1}{6}\ \ &((i,j)\in\mathcal I_6),
\\
0\ \ &((i,j)\notin\mathcal I_6),
\end{cases}
\ \ (1\leq i<j\leq k),
\end{aligned}
\end{align}
where
\begin{align*}
\mathcal I_6=\{(i,i+1):1\leq i\leq5\}\cup\{(1,6)\}.
\end{align*}

\textbf{Step 2. Absolute proximity to a regular hexagon.}

By \eqref{abomegabubun} and \eqref{akbkkeisuu}, we have for $1\leq i\leq 6$ and $(i,j)\in\mathcal{I}_6$ 
\begin{align*}
\lim_{n\to\infty}\frac{\Psi(\rho_i)}{\mathcal{A}}=\lim_{n\to\infty}\frac{\Psi(\rho_{ij})}{\mathcal{A}}=\frac{1}{6}.
\end{align*}
In addition, by the definition of $D$ and \eqref{Psies0922}, we obtain
\begin{align}\label{nagasaes1005}
\begin{aligned}
\rho_i&=D+o(1)\ \ (1\leq i\leq6),
\\
\rho_{ij}&=D+o(1)\ \ ((i,j)\in\mathcal I_6),
\\
\rho_{il}-D&\to\infty\ \  (1\leq i\leq6,\ 7\leq l\leq k).
\end{aligned}
\end{align}
Furthermore, by the law of cosines, we have for $(i,j)\in\mathcal{I}_6$
\begin{align*}
\frac{(z_i-z_0)\cdot(z_j-z_0)}{\rho_i\rho_j}=\frac{\rho_i^2+\rho_j^2-\rho_{ij}^2}{2\rho_i\rho_j}=\frac{1}{2}+o(D^{-1}).
\end{align*}
Together with the cyclic ordering of the limiting hexagon, this shows that the successive angular gaps are $\frac{\pi}{3}+o(D^{-1})$. Thus, we have 
\begin{align}\label{rho13order}
\rho_{13}=\sqrt3D+o(1).
\end{align}
Furthermore, there exist rotations $R_n\in SO(2)$ such that
\begin{align*}
\lim_{n\to\infty}\max_{1\leq i\leq6}|z_i^{(n)}-z_0^{(n)}-D_nR_nu_i|=0.
\end{align*}
Here, let $q$ be the minimizing hexagonal configuration
associated with $(z_0,\ldots,z_6)$, and set
\begin{align*}
w=(z_0,\ldots,z_6)-q,\ r=|q_1-q_0|.
\end{align*}
By minimality of $q$, the preceding proximity estimate gives
\begin{align}\label{Lemma4.4tekiyou}
|w|=o(1),\qquad r=D+o(1),\ g(r)\sim\mathcal A.
\end{align}

\textbf{Step 3. The remaining interactions.}

Next, we show that we have for $7\leq j\leq k$
\begin{align*}
\lim_{n\to\infty}\left(\rho_j-\sqrt{3}D\right)=\infty.
\end{align*}
Fix $j$ with $7\leq j\leq k$ and set
\begin{align*}
m_j=\min_{1\leq i\leq6}|z_j-z_i|.
\end{align*}
By \eqref{nagasaes1005}, we have $\lim_{n\to\infty}(m_j-D)=\infty$. We choose a vertex direction $R_nu_i$ whose angle with $z_j-z_0$ is at most $\frac{\pi}{6}$. Then, we have
\begin{align}\label{seirokusotoes1}
m_j\leq \sqrt{\rho_j^2+D^2-\sqrt3D\rho_j}+o(1).
\end{align}
Here, we fix $M>0$ and we assume that there exists a sequence $(\tilde{z}^{(n)})_n\subset (z^{(n)})_n$ such that 
\begin{align*}
\lim_{n\to\infty}\left(m_j(\tilde{z}^{(n)})-D(\tilde{z}^{(n)})\right)=\infty,\ \ \rho_j(\tilde{z}^{(n)})\leq \sqrt{3}D(\tilde{z}^{(n)})+M.
\end{align*}
Then, since we have $\rho_j\geq D$, by \eqref{seirokusotoes1} we obtain
\begin{align*}
m_j\leq\sqrt{D^2+\sqrt3DM+M^2}+o(1)=D+O_M(1),
\end{align*}
which contradicts $m_j-D\to\infty$. Therefore, we have $\lim_{n\to\infty}(\rho_j-\sqrt{3}D)=\infty$, and hence we obtain
\begin{align}\label{seirokusotoes2}
\sum_{j=7}^k\Psi(\rho_j)=o\left(\Psi(\sqrt{3}D)\right).
\end{align}
We note that when we have $k=6$, this sum is understood to be zero. Here, we define 
\begin{align}\label{tildeH0Pdef}
\begin{aligned}
&\tilde{H}_0(z_0,z_1,\ldots,z_k)=H_0(z_0,z_1,\ldots,z_6),
\\
&P(z)=\sum_{1\leq i<j\leq k,\ (i,j)\notin\mathcal{I}_6}\Psi(\rho_{ij})+\sum_{i=1}^k\mathcal N_{-1}(\rho_i) +\sum_{1\leq i<j\leq k}\mathcal N_1(\rho_{ij})-\sum_{l=7}^k\Psi(\rho_l).
\end{aligned}
\end{align}
Then, by \eqref{Gdef}, we have 
\begin{align}\label{G=H-P}
G=\tilde{H}_0-P.
\end{align}
Furthermore, by \eqref{mathcalNsei}, \eqref{rho13order}, and \eqref{seirokusotoes2}, we obtain
\begin{align}\label{Pgtrsim}
P\geq\Psi(\rho_{13})-\sum_{j=7}^k\Psi(\rho_j)=(1+o(1))\Psi(\sqrt{3}D)\gtrsim \Psi(\sqrt{3}D).
\end{align}
Furthermore, by \eqref{Psies0922} and \eqref{mathcalNbibunes}, there exists $C_P>0$ such that 
\begin{align}\label{Plesssim}
|\nabla_zP|\leq C_PP.
\end{align}

\textbf{Step 4. Application of Lemma \ref{H0nogosalemma}.}

Finally, we derive a contradiction. By \eqref{betastardef}, \eqref{prop4.1nokatei}, \eqref{G=H-P}, and \eqref{Pgtrsim}, we have for all sufficiently large $n$,
\begin{align}\label{tildeHtoPes}
\tilde{H}_0=G+P\geq P-ce^{-\beta_\star D}\geq\frac{1}{2}P>0.
\end{align}
By \eqref{Apsigsim}, \eqref{Lemma4.4tekiyou}, and Lemma \ref{H0nogosalemma}, there exists $C_{\tilde{H}_0}>0$ such that
\begin{align}\label{Lemma4.4tukatta}
\begin{aligned}
\tilde{H}_0&\leq C_{\tilde{H}_0}\mathcal{A}|w|^2,
\\
|\nabla_z\tilde{H}_0|&\geq \frac{1}{C_{\tilde{H}_0}}\mathcal{A}|w|.
\end{aligned}
\end{align}
By \eqref{Lemma4.4tekiyou}, \eqref{Plesssim}, \eqref{tildeHtoPes}, and \eqref{Lemma4.4tukatta}, we obtain
\begin{align}\label{nablaGlower1}
|\nabla_zG|&\geq|\nabla_z\tilde{H}_0|-|\nabla_zP|\geq \frac{1}{C_{\tilde{H}_0}}\mathcal{A}|w|-2C_PC_{\tilde{H}_0}\mathcal{A}|w|^2\gtrsim\mathcal A|w|.
\end{align}
Furthermore, by \eqref{G=H-P}, \eqref{Pgtrsim}, \eqref{tildeHtoPes}, and \eqref{Lemma4.4tukatta}, we obtain
\begin{align}\label{nablaGlower2}
\max{(G,0)}+\Psi(\sqrt3D)\lesssim \tilde{H}_0\lesssim\mathcal A|w|^2.
\end{align}
By \eqref{nablaGlower1} and \eqref{nablaGlower2}, we obtain
\begin{align*}
|\nabla_zG|^2\gtrsim\mathcal A^2|w|^2\gtrsim\mathcal A\left(\max(G,0)+\Psi(\sqrt3D)\right),
\end{align*}
which contradicts \eqref{prop4.1nokatei}. Therefore, we complete the proof.
\end{proof}

\section{The Cauchy--Schwarz deficit and constraints on soliton configurations}

In this section, we use the argument based on the deficit in the Cauchy--Schwarz inequality developed in \cite{I4} to derive the limiting balance equations along a sequence of times. We then prove Theorem \ref{maintheorem} using these balance equations.

\subsection{Positivity of the interaction energy}

In this subsection, we show that $\mathcal{U}$ is comparable to $g(D)$, which allows us to apply the argument based on the deficit in the Cauchy--Schwarz inequality.
\begin{lemma}\label{energysim}
Let $\vec{u}$ be a $(1,k)$-sign soliton solution satisfying \eqref{1,ksolicondrevise}. Then, we have for $t\gg 1$
\begin{align}\label{GsimUsimg}
G(t)\sim g(D(t)),\ \ \left|\nabla G(t)\right|\sim g(D(t)).
\end{align}
\end{lemma}
\begin{proof}
Throughout the proof, all quantities are evaluated at $z(t)$, and $t$ is sufficiently large.

By \eqref{Psies0922}, \eqref{Genergylower}, and Lemma \ref{mathcalNeslem}, we have
\begin{align}\label{AsimPsisimgsec5}
\begin{aligned}
&\mathcal A\sim\Psi(D)\sim g(D),
\\
&|\nabla_z\mathcal A|+|\nabla_zG|\lesssim \mathcal{A}.
\end{aligned}
\end{align}
First, we prove $G(t)>0$. By Proposition \ref{haitipropo} and \eqref{betastardef}, we have
\begin{align}\label{Prop41katuyou}
e^{-\frac{\beta_\star+\sqrt{3}}{2} D}\lesssim \sqrt{\mathcal A\Psi(\sqrt3D)}\lesssim |\nabla_zG|.
\end{align}
By \eqref{zesref0927}  we have
\begin{align}\label{Gnobibun1008}
\frac{d}{dt}G=\nabla_zG\cdot\dot z=-\frac{|\nabla_zG|^2}{C_{\star}}+O\left(e^{-\beta_\star D}|\nabla_zG|\right).
\end{align}
By  \eqref{Prop41katuyou} and \eqref{Gnobibun1008}, we obtain
\begin{align}\label{Gtozbibun}
\frac{d}{dt}G\sim-|\nabla_zG|^2,\ \ |\dot z|\lesssim|\nabla_zG|.
\end{align}
Therefore, $G$ is decreasing. Furthermore, since we have $\lim_{t\to\infty}G(t)=0$, we have $G(t)>0$. In addition,  by Proposition \ref{haitipropo}, we obtain
\begin{align}\label{Prop41kairyou}
\sqrt{\mathcal{A}G}\lesssim |\nabla_zG|.
\end{align}

Next, we prove $\mathcal{A}\lesssim G$. By \eqref{AsimPsisimgsec5} and \eqref{Gtozbibun}, we have
\begin{align}\label{Airoiro1008}
\left|\frac{d}{dt}\mathcal{A}\right|\leq|\nabla_z\mathcal{A}||\dot z|\lesssim\mathcal{A}|\nabla_zG|.
\end{align}
We define
\begin{align*}
h(t)=\sqrt{\frac{G(t)}{\mathcal{A}(t)}}.
\end{align*}
Then, by \eqref{Gtozbibun}, \eqref{Prop41kairyou}, and \eqref{Airoiro1008}, there exists $0<c\ll 1$ such that
\begin{align}\label{hnobs}
\dot{h}=\frac{\dot{G}}{2\sqrt{\mathcal{A}G}}-\frac{h\dot{\mathcal{A}}}{2\mathcal{A}}\leq\left(-c+\frac{h}{c}\right)|\nabla_zG|.
\end{align}
Here, we introduce the bootstrap estimate
\begin{align}\label{sec5bs}
h(t)\leq \frac{c^2}{2}.
\end{align}
We assume that  there exists $T\gg 1$ such that 
\begin{align}\label{bskatei}
h(T)\leq \frac{c^2}{2}.
\end{align}
Then, we define $T_1\in [T,\infty]$ as 
\begin{align}\label{T1def}
T_1=\sup{ \{ t\in [T,\infty)\ \mbox{such that}\ \eqref{sec5bs}\ \mbox{holds on}\ [T,t]\}}.
\end{align}
Furthermore, we assume $T_1<\infty$. Then, by \eqref{hnobs}, we have $\dot{h}(T_1)<0$, which contradicts the maximality of $T_1$. Thus, we obtain $T_1=\infty$. Then, we have for $t>T$
\begin{align}\label{dothkairyou}
\dot{h}\leq -\frac{c}{2}|\nabla_zG|.
\end{align}
Then, by \eqref{Gtozbibun} and \eqref{dothkairyou}, we have for $T\ll t$,
\begin{align*}
D(t)-D(T)\lesssim \int_{T}^{t}|\dot z(s)|ds\lesssim\int_{T}^{t}|\nabla_zG(z(s))|ds\lesssim h(T),
\end{align*}
which contradicts $D(t)\to\infty$. Thus, we have $h(t)\geq \frac{c^2}{2}$. Furthermore, since we have $G\leq \mathcal{A}$, by \eqref{AsimPsisimgsec5} we have $G\sim\mathcal{A}\sim g(D)$. In addition, $|\nabla_zG|\sim \mathcal{A}$ follows from \eqref{Gdef} and \eqref{Prop41kairyou}, and therefore we obtain \eqref{GsimUsimg}.
\end{proof}

\subsection{Application of the Cauchy--Schwarz deficit argument}

In this subsection, we apply \cite[Section 4]{I4}. First, we define
\begin{align}
\mathcal{S}&=\frac{1}{C_\star}\left|\nabla G\right|^2, \label{mathcalSdef}
\\
\mathcal{L}&=-\log{G},\label{mathcalLdef}
\\
z_g&=\frac{1}{k+1}\sum_{i=0}^kz_i,  \label{zgdef}
\\
\mathcal{I}&=C_{\star}\sum_{i=0}^k|z_i-z_g|^2, \label{mathcalIdef}
\\
\mathcal{M}&=-\sum_{i=0}^k(z_i-z_g)\cdot \nabla_{z_i}G. \label{mathcalMdef}
\end{align}
Then, by Lemma \ref{energysim}, we have
\begin{align}\label{L-Dsim}
|\mathcal{L}(t)-D(t)|&\lesssim \log{D(t)}.
\end{align}
Furthermore, we have the following lemma.
\begin{lemma}\label{DLIodenado}
Let $\vec{u}$ be a $(1,k)$-sign soliton solution satisfying \eqref{1,ksolicondrevise}. Then, we have for $t\gg 1$
\begin{align}
\dot{G}&=-\mathcal{S}+O\left(\frac{G^2}{D}\right), \label{dotGsim}
\\
\dot{\mathcal{L}}&=\frac{\mathcal{S}}{G}+O\left(\frac{G}{D}\right), \label{dotLsim}
\\
\frac{1}{2}\dot{\mathcal{I}}&=\mathcal{M}+O\left(\frac{\mathcal{I}^{\frac{1}{2}}G}{D}\right), \label{dotIsim}
\\
\mathcal{M}&=DG+O(G). \label{M=DGkurai}
\end{align}

\end{lemma}

\begin{proof}
We prove these estimates using essentially the same computations as in \cite[Lemmas 4.1 and 4.2]{I4}. We therefore refer the reader to \cite{I4} for the detailed computations. By Lemma \ref{chushinseimitulem} and Lemma \ref{energysim}, we have
\begin{align}\label{chusinodesec5}
\dot{z}_i=-\frac{1}{C_{\star}}\nabla_{z_i}G+O\left(\frac{G}{D}\right).
\end{align}
By \eqref{GsimUsimg} and \eqref{chusinodesec5}, we have
\begin{align*}
\dot{G}&=\sum_{i=0}^k\nabla_{z_i}G\cdot \dot{z}_i=-\frac{1}{C_{\star}}|\nabla_zG|^2+O\left(\frac{G|\nabla_zG|}{D}\right)=-\mathcal{S}+O\left(\frac{G^2}{D}\right),
\end{align*}
which implies \eqref{dotGsim}. In addition, by \eqref{dotGsim}, we have
\begin{align*}
\dot{\mathcal{L}}=-\frac{\dot{G}}{G}=\frac{\mathcal{S}}{G}+O\left(\frac{G}{D}\right),
\end{align*}
and we obtain \eqref{dotLsim}. Next, we estimate \eqref{dotIsim}. Before we estimate \eqref{dotIsim}, we note that we have $\sum_{i=0}^k(z_i-z_g)=0$. Then, by direct computation, we have
\begin{align}\label{keisan1002-1}
\begin{aligned}
\frac{1}{2}\dot{\mathcal{I}}&=C_{\star}\sum_{i=0}^k(z_i-z_g)\cdot (\dot{z}_i-\dot{z}_g)
\\
&=C_{\star}\sum_{i=0}^k(z_i-z_g)\cdot \dot{z}_i
\\
&=\mathcal{M}+O\left(\frac{G}{D}\sum_{i=0}^k|z_i-z_g|\right).
\end{aligned}
\end{align}
Furthermore, by the Cauchy--Schwarz inequality, we have
\begin{align}\label{keisan1002-2}
\sum_{i=0}^k|z_i-z_g|\leq \sqrt{k+1}\left(\sum_{i=0}^k|z_i-z_g|^2\right)^{\frac{1}{2}}.
\end{align}
By \eqref{keisan1002-1} and \eqref{keisan1002-2}, we obtain \eqref{dotIsim}. Last, we estimate \eqref{M=DGkurai}. By direct computation, using \eqref{Gdef} and \eqref{mathcalMdef}, we have
\begin{align*}
\mathcal{M}&=\sum_{i=1}^k\rho_i g(\rho_i)-\sum_{1\leq i<j\leq k}\rho_{ij}g(\rho_{ij})+\sum_{i=1}^k\rho_i\mathcal N_{-1}^{\prime}(\rho_i)+\sum_{1\leq i<j\leq k}\rho_{ij}\mathcal N_1^{\prime}(\rho_{ij}).
\end{align*}
Therefore, we obtain
\begin{align}\label{1002keisan3}
\begin{aligned}
\mathcal{M}-DG&=\sum_{i=1}^k\left(\rho_ig(\rho_i)-D\Psi(\rho_i)\right)-\sum_{1\leq i<j\leq k}\left(\rho_{ij}g(\rho_{ij})-D\Psi(\rho_{ij})\right)
\\
&\quad+\sum_{i=1}^k\left(\rho_i\mathcal{N}_{-1}^{\prime}(\rho_i)+D\mathcal{N}_{-1}(\rho_i)\right)+\sum_{1\leq i<j\leq k}\left(\rho_{ij}\mathcal{N}_1^{\prime}(\rho_{ij})+D\mathcal{N}_1(\rho_{ij})\right).
\end{aligned}
\end{align}
Furthermore, we have for $r\geq D\gg 1$
\begin{align}\label{1002keisan4}
\begin{aligned}
\left|rg(r)-D\Psi(r)\right|&\leq (r-D)g(r)+D|g(r)-\Psi(r)|
\\
&\lesssim (r-D+1)e^{-(r-D)}g(D)\lesssim g(D),
\end{aligned}
\end{align}
and by Lemma \ref{mathcalNeslem}, we have for $s=\pm 1$ and $r\geq D\gg 1$
\begin{align}\label{1002keisan5}
|r\mathcal{N}_s^{\prime}(r)|+|D\mathcal{N}_s(r)|\lesssim (r+D)e^{-\gamma r}\lesssim g(D),
\end{align}
where $\gamma>1$ is a constant. By \eqref{GsimUsimg}, \eqref{1002keisan3}, \eqref{1002keisan4}, and \eqref{1002keisan5}, we obtain \eqref{M=DGkurai}, and we complete the proof.
\end{proof}

Using Lemma \ref{DLIodenado}, we also obtain the following estimates.

\begin{lemma}\label{SsimIsimlem}
Let $\vec{u}$ be a $(1,k)$-sign soliton solution satisfying \eqref{1,ksolicondrevise}. Then, we have for $t\gg 1$
\begin{align}
\mathcal{S}\sim G^2, \label{SsimG}
\\
\mathcal{I}\sim \mathcal{L}^2. \label{IsimL^2}
\end{align}
\end{lemma}
\begin{proof}
By Lemma \ref{energysim} and \eqref{mathcalSdef}, we have \eqref{SsimG}. By \eqref{dotLsim} and \eqref{SsimG}, we have
\begin{align}\label{dotLsimkai}
\dot{\mathcal{L}}\sim G.
\end{align}
By the Cauchy--Schwarz inequality, we have
\begin{align}\label{ko-shi-}
\mathcal{M}^2\leq \mathcal{I}\mathcal{S}.
\end{align}
By \eqref{dotIsim}, \eqref{dotLsimkai}, and \eqref{ko-shi-}, we have
\begin{align*}
\left|\frac{d}{dt}\sqrt{\mathcal{I}}\right|\lesssim \frac{|\mathcal{M}|}{\sqrt{\mathcal{I}}}+\frac{G}{D}\lesssim \sqrt{\mathcal{S}}+\frac{G}{D}\sim G\sim \dot{\mathcal{L}},
\end{align*}
and we obtain
\begin{align}\label{Ilesssimhyouka}
\mathcal{I}\lesssim \mathcal{L}^2.
\end{align}
Furthermore, since we have for $i\neq j$,
\begin{align*}
D^2\leq |z_i-z_j|^2\leq 2|z_i-z_g|^2+2|z_j-z_g|^2,
\end{align*}
by \eqref{L-Dsim}, we obtain
\begin{align}\label{Igtrsim}
\mathcal{L}^2\lesssim \mathcal{I}.
\end{align}
By \eqref{Ilesssimhyouka} and \eqref{Igtrsim}, we obtain \eqref{IsimL^2}.
\end{proof}

Using the preceding estimates, we show that the normalized Cauchy--Schwarz deficit tends to zero along a sequence of times.
\begin{lemma}\label{kessinlem}
Let $\vec{u}$ be a $(1,k)$-sign soliton solution satisfying \eqref{1,ksolicondrevise}. Then, there exists a sequence of times $t_n\to\infty$ such that 
\begin{align}\label{kessonkyoku}
\lim_{n\to\infty}\frac{\mathcal{I}(t_n)\mathcal{S}(t_n)-\mathcal{M}^2(t_n)}{\mathcal{L}^2(t_n)G^2(t_n)}=0.
\end{align}
\end{lemma}
\begin{proof}
This lemma also follows from essentially the same direct computations as in \cite[Lemma 4.3]{I4}. We therefore refer the reader to \cite{I4} for the details. We define
\begin{align*}
\mathcal{P}=\frac{\mathcal{I}}{\mathcal{L}^2},\ \mathcal{Q}=\frac{\mathcal{S}}{G^2},\ \mathcal{O}=\frac{\mathcal{M}}{\mathcal{L}G}.
\end{align*}
Then, we have
\begin{align}\label{PQOsim}
\mathcal{P}\sim 1,\ \mathcal{Q}\sim 1,\ |\mathcal{O}-1|\lesssim \frac{\log{\mathcal{L}}}{\mathcal{L}}.
\end{align}
By direct computation, we also have
\begin{align}\label{PtoLsim}
\frac{\mathcal{L}}{\dot{\mathcal{L}}}\dot{\mathcal{P}}=\frac{2\mathcal{O}}{\mathcal{Q}}-2\mathcal{P}+O\left(\frac{1}{\mathcal{L}}\right).
\end{align}
Here, we define $\mathcal{J}=\mathcal{P}\mathcal{Q}-\mathcal{O}^2$. Then, by \eqref{ko-shi-}, we have $\mathcal{J}\geq 0$. Furthermore, by \eqref{PQOsim} and \eqref{PtoLsim}, we have
\begin{align*}
\frac{\mathcal{L}}{\dot{\mathcal{L}}}\dot{\mathcal{P}}=-\frac{2\mathcal{J}}{\mathcal{Q}}+O\left(\frac{\log{\mathcal{L}}}{\mathcal{L}}\right).
\end{align*}
Thus, there exists $c>0$ such that, for all sufficiently large $t$,
\begin{align}\label{Pbibunhutoushiki}
\dot{\mathcal{P}}\leq -c\mathcal{J}\frac{\dot{\mathcal{L}}}{\mathcal{L}}+\frac{1}{c}\frac{\log{\mathcal{L}}}{\mathcal{L}^2}\dot{\mathcal{L}}.
\end{align}
Integrating \eqref{Pbibunhutoushiki} on $[t,T]$, where $1\ll t<T$, we have
\begin{align*}
\mathcal{P}(T)-\mathcal{P}(t)&\leq -c\int_t^T\mathcal{J}(s)\frac{\dot{\mathcal{L}}(s)ds}{\mathcal{L}(s)}
\\
&\quad+\frac{1}{c}\left\{\left(\frac{\log{\mathcal{L}(t)}}{\mathcal{L}(t)}+\frac{1}{\mathcal{L}(t)}\right)-\left(\frac{\log{\mathcal{L}(T)}}{\mathcal{L}(T)}+\frac{1}{\mathcal{L}(T)}\right)\right\},
\end{align*}
and we obtain
\begin{align}\label{JtoLsekibunyuukai}
\int_t^{\infty}\mathcal{J}(s)\frac{\dot{\mathcal{L}}(s)ds}{\mathcal{L}(s)}<\infty.
\end{align}
By \eqref{JtoLsekibunyuukai} and $\int_t^{\infty}\frac{\dot{\mathcal{L}}(s)ds}{\mathcal{L}(s)}=\infty$, there exists a sequence $t_n\to\infty$ such that 
\begin{align*}
\lim_{n\to\infty}\mathcal{J}(t_n)=0,
\end{align*}
which implies \eqref{kessonkyoku}.
\end{proof}

\subsection{Constraints on the limiting configuration}

In this subsection, we pass to a further subsequence of the times given by Lemma \ref{kessinlem} to derive the limiting balance equations.
\begin{lemma}\label{turiaieq1002}
Let $\vec{u}$ be a $(1,k)$-sign soliton solution satisfying \eqref{1,ksolicondrevise}. Then, there exist $t_n\to\infty$, $y_1,y_2,\ldots, y_k\in \mathbb{R}^2$ and $c>0$ and $\alpha_i\geq 0$ ($1\leq i\leq k$) and $\beta_{ij}\geq 0$ ($1\leq i\neq j\leq k$) such that  
\begin{align}\label{tenretu1002}
\begin{aligned}
y_i&=\lim_{n\to\infty}\frac{z_i(t_n)-z_0(t_n)}{\mathcal{L}(t_n)},
\\
\frac{1}{c}&=\lim_{n\to\infty}\frac{\mathcal{I}(t_n)}{C_{\star}\mathcal{L}(t_n)^2},
\\
\alpha_i&=\lim_{n\to\infty}\frac{g(\rho_i(t_n))}{G(t_n)},
\\
\beta_{ij}&=\lim_{n\to\infty}\frac{g(\rho_{ij}(t_n))}{G(t_n)}.
\end{aligned}
\end{align}
Furthermore, we have
\begin{align}\label{haichi1002-1}
\begin{aligned}
|y_1|=|y_2|=\cdots=|y_k|=1,\ |y_i-y_j|\geq 1\ (i\neq j),
\\
\beta_{ij}=\beta_{ji},
\\
\beta_{ij}=0\ \ (\mbox{if}\ |y_i-y_j|>1),
\\
cy_i=\alpha_iy_i+\sum_{j=1}^k\alpha_jy_j+\sum_{j\neq i}\beta_{ij}(y_j-y_i)\ \ (1\leq i\leq k).
\end{aligned}
\end{align}
\end{lemma}
\begin{proof}
Let $\{t_n\}$ be the sequence given by Lemma \ref{kessinlem}. First, by Lemma \ref{energysim}, \eqref{mathcalIdef}, and \eqref{IsimL^2}, we have
\begin{align}\label{simdeyuukairetu}
\begin{aligned}
0<\frac{|z_i-z_0|}{\mathcal{L}}&\lesssim \frac{\sqrt{\mathcal{I}}}{\mathcal{L}}\lesssim 1,
\\
0<\frac{g(\rho_i)}{G}&\lesssim \frac{g(D)}{G}\lesssim 1,
\\
0< \frac{g(\rho_{ij})}{G}&\lesssim \frac{g(D)}{G}\lesssim 1.
\end{aligned}
\end{align}
By \eqref{simdeyuukairetu}, passing to a subsequence and keeping the same notation, we obtain \eqref{tenretu1002}, with  $y_i$, $c>0$, $\alpha_i\geq0$, and $\beta_{ij}\geq0$. Since we have $\rho_{ij}=\rho_{ji}$, we also have $\beta_{ij}=\beta_{ji}$. Furthermore, by \eqref{L-Dsim}, we obtain
\begin{align}\label{coefficientsupport1002}
\alpha_i=0\ \  (|y_i|>1),\ \beta_{ij}=0\ \  (|y_i-y_j|>1).
\end{align}

Next, we derive the limiting balance equation. By \eqref{mathcalSdef}, \eqref{mathcalIdef}, \eqref{mathcalMdef}, \eqref{IsimL^2}, and \eqref{kessonkyoku}, we have
\begin{align}\label{deficitalignment1002}
\lim_{n\to\infty}\sum_{i=0}^k\left|\frac{\nabla_{z_i}G}{G}+\frac{C_\star\mathcal{M}}{\mathcal{I}G}(z_i-z_g)\right|^2=\lim_{n\to\infty}\frac{C_\star}{\mathcal{I}G^2}\left(\mathcal{I}\mathcal{S}-\mathcal{M}^2\right)=0.
\end{align}
Moreover, by \eqref{L-Dsim} and \eqref{M=DGkurai}, we have
\begin{align}\label{hiritu1010}
\lim_{n\to\infty}\frac{\mathcal{M}}{\mathcal{L}G}=1,\ \ \lim_{n\to\infty}\frac{C_\star\mathcal{M}\mathcal{L}}{\mathcal{I}G}
=c.
\end{align}
By \eqref{deficitalignment1002} and \eqref{hiritu1010}, and
\begin{align*}
&\quad -\frac{\nabla_{z_i}G-\nabla_{z_0}G}{G}
\\
&=\left( \frac{\nabla_{z_0}G}{G}+\frac{C_{\star}\mathcal{M}}{\mathcal{I}G}(z_0-z_g)\right)-\left(\frac{\nabla_{z_i}G}{G}+\frac{C_\star\mathcal{M}}{\mathcal{I}G}(z_i-z_g)\right)+\frac{C_{\star}\mathcal{M}}{\mathcal{I}G}Z_i,
\end{align*}
we obtain
\begin{align}\label{gradientlimit1002}
cy_i=\lim_{n\to\infty}-\frac{\nabla_{z_i}G(t_n)-\nabla_{z_0}G(t_n)}{G(t_n)}.
\end{align}
In addition, by Lemma \ref{mathcalNeslem}, we have
\begin{align*}
\left|\nabla_z\mathcal{N}_{-1}(\rho_i)\right|+\left|\nabla_z\mathcal{N}_1(\rho_{ij})\right|=o(G),
\end{align*}
and therefore we obtain
\begin{align}\label{ynokyokugensiki}
-\frac{\nabla_{z_i}G-\nabla_{z_0}G}{G}=\frac{g(\rho_i)}{G}\frac{Z_i}{\rho_i}+\sum_{j=1}^k\frac{g(\rho_j)}{G}\frac{Z_j}{\rho_j}+
\sum_{\substack{1\leq j\leq k,\ j\neq i}}\frac{g(\rho_{ij})}{G}\frac{Z_{ij}}{\rho_{ij}}+o(1).
\end{align}
By \eqref{gradientlimit1002} and \eqref{ynokyokugensiki}, we obtain
\begin{align}\label{cyrutiaikari}
cy_i=\alpha_i\frac{y_i}{|y_i|}+\sum_{j=1}^k\alpha_j\frac{y_j}{|y_j|}+\sum_{1\leq j\leq k,\ j\neq i}\beta_{ij}\frac{y_j-y_i}{|y_j-y_i|}.
\end{align}
By \eqref{coefficientsupport1002} and \eqref{cyrutiaikari}, we obtain
\begin{align}\label{limitingbalance1002proof}
cy_i=\alpha_i y_i+\sum_{j=1}^k\alpha_j y_j+\sum_{j\neq i}\beta_{ij}(y_j-y_i).
\end{align}

Finally, we prove that $|y_i|=1$ for every $i$. Since it is clear that $|y_i|\geq1$ for every $i$, it suffices to show that $\max_{1\leq i\leq k}|y_i|\leq1$. We assume that we have
\begin{align*}
R:=\max_{1\leq i\leq k}|y_i|>1.
\end{align*}
Then, we choose $i_0$ such that $|y_{i_0}|=R$, and we define $u=\frac{y_{i_0}}{R}$. Since we have $\alpha_{i_0}=0$, taking the inner product of \eqref{limitingbalance1002proof} for $i=i_0$ with $u$, we obtain
\begin{align}\label{cRtukuru}
cR=\sum_{j=1}^k\alpha_j y_j\cdot u+\sum_{j\neq i_0}\beta_{i_0j}(y_j\cdot u-R)\leq \sum_{j=1}^k\alpha_j y_j\cdot u.
\end{align}
By \eqref{limitingbalance1002proof}, we have
\begin{align}\label{ynowaeq}
c\sum_{i=1}^k y_i=(k+1)\sum_{i=1}^k\alpha_iy_i.
\end{align}
By \eqref{cRtukuru} and \eqref{ynowaeq}, we obtain
\begin{align*}
cR\leq \sum_{j=1}^k\alpha_j y_j\cdot u=\frac{c}{k+1}\sum_{j=1}^ky_j\cdot u\leq \frac{ckR}{k+1}<cR,
\end{align*}
which is a contradiction. Therefore, we have $R=1$ and we have for every $i$ $|y_i|=1$. Thus, we obtain \eqref{haichi1002-1}.
\end{proof}

We now prove Theorem \ref{maintheorem} using Lemma \ref{turiaieq1002}.
\begin{proof}[Proof of Theorem \ref{maintheorem}]
By the sign-change invariance in \eqref{invariance}, it suffices to prove the nonexistence of $(1,k)$-sign soliton solutions.
Suppose, for a contradiction, that $\vec{u}$ is a $(1,k)$-sign soliton solution for some $k\geq6$. By suitably choosing the centers, we may assume that \eqref{1,ksolicondrevise} holds. We distinguish the following two cases.

\textbf{Case 1: $k\geq7$.}

By Lemma \ref{turiaieq1002}, there exist $y_1,y_2,\ldots, y_k\in S^1$ satisfying \eqref{haichi1002-1}. Then, since we have $y_1,y_2,\ldots, y_k\in S^1$, we obtain
\begin{align*}
\min_{1\leq i<j\leq k}|y_i-y_j|\leq 2\sin{\frac{\pi}{k}}\leq 2\sin{\frac{\pi}{7}}<1,
\end{align*}
which contradicts \eqref{haichi1002-1}. Therefore, no $(1,k)$-sign soliton solution exists for $k\geq7$.

\textbf{Case 2: $k=6$.}

By Lemma \ref{turiaieq1002}, there exist $t_n$, $y_i$, $\alpha_i$, $\beta_{ij}$, and $c>0$ satisfying \eqref{tenretu1002} and \eqref{haichi1002-1}. In particular, by \eqref{haichi1002-1} and the condition $|y_i-y_j|\geq1$, after suitably relabeling the outer centers $z_1,\ldots,z_6$ and applying a rotation, we may assume that
\begin{align}\label{yrokkaku}
y_i=\left(\cos\frac{(i-1)\pi}{3},\sin\frac{(i-1)\pi}{3}\right)\ \  (1\leq i\leq6).
\end{align}
In the remainder of this proof, the indices of the outer centers $z_1,\ldots,z_6$ are understood modulo $6$. We define $\beta_i=\beta_{i,i+1}$ for $i\in\mathbb{Z}/6\mathbb{Z}$. Then, by \eqref{haichi1002-1}, \eqref{ynowaeq}, and \eqref{yrokkaku}, we have for $i\in\mathbb{Z}/6\mathbb{Z}$
\begin{align}\label{turiairokkakukei}
cy_i=\alpha_i y_i+\beta_i(y_{i+1}-y_i)+\beta_{i-1}(y_{i-1}-y_i).
\end{align}
Since we have $y_{i-1}+y_{i+1}=y_i$ for $i\in \mathbb{Z}/6\mathbb{Z}$, by \eqref{turiairokkakukei}, we have for $i\in\mathbb{Z}/6\mathbb{Z}$,
\begin{align}\label{yi-1yi+1itijidokuritu}
(c-\alpha_i)(y_{i-1}+y_{i+1})=(c-\alpha_i)y_i=-\beta_iy_{i-1}-\beta_{i-1}y_{i+1}.
\end{align}
Since $y_{i-1}$ and $y_{i+1}$ are linearly independent, by \eqref{yi-1yi+1itijidokuritu}, we have for $i\in\mathbb{Z}/6\mathbb{Z}$
\begin{align*}
\alpha_i-\beta_i=c,\ \  \beta_i=\beta_{i-1}.
\end{align*}
Therefore, there exist $a>0$ and $b\geq 0$ such that for every $i\in\mathbb{Z}/6\mathbb{Z}$
\begin{align}\label{hexagonalcoefficientineq}
\alpha_i=a,\ \ \beta_i=b,\ \ a=b+c>b.
\end{align}
Then, by \eqref{tenretu1002} and \eqref{hexagonalcoefficientineq}, we have for $i\in\mathbb{Z}/6\mathbb{Z}$
\begin{align}\label{rhoinnagasahabtei1}
\lim_{n\to\infty}\frac{g(\rho_{i,i+1}(t_n))}{G(t_n)}=b,\ \ \lim_{n\to\infty}\frac{g(\rho_i(t_n))}{G(t_n)}=a.
\end{align}
In particular, by \eqref{gnozenkin1}, \eqref{hexagonalcoefficientineq}, and \eqref{rhoinnagasahabtei1}, we have for $n\gg 1$ and $i\in\mathbb{Z}/6\mathbb{Z}$
\begin{align}\label{hexagonedgelower}
\rho_{i,i+1}(t_n)>\max\left\{\rho_i(t_n),\rho_{i+1}(t_n)\right\}.
\end{align}
Here, we fix a sufficiently large $n$ and define $\theta_i$ to be the angle at $z_0$ in the triangle with vertices $z_0,z_i,z_{i+1}$, with all centers evaluated at $t_n$. We choose $\tilde{i}$ such that $\theta_{\tilde{i}}=\min_{1\leq i\leq6}\theta_i$. Since we have $\sum_{i=1}^6\theta_i=2\pi$, we obtain $\theta_{\tilde{i}}\leq\frac{\pi}{3}$. Therefore, in the triangle with vertices $z_0,z_{\tilde{i}},z_{\tilde{i}+1}$, the angle at $z_0$ is not the unique largest angle, and hence we have
\begin{align*}
\rho_{\tilde{i},\tilde{i}+1}^2\leq \max\{\rho_{\tilde{i}},\rho_{\tilde{i}+1}\}^2,
\end{align*}
which contradicts \eqref{hexagonedgelower}. Therefore, no $(1,k)$-sign soliton solution exists in the case $k=6$ either.

By Case 1 and Case 2, we complete the proof of Theorem \ref{maintheorem}.
\end{proof}

\section*{Statements and Declarations}

\textbf{Funding} No funds, grants, or other support was received.

\textbf{Competing interests} The author has no relevant financial or non-financial interests to disclose.

\textbf{Data availability} No datasets were generated or analyzed during the current study.

\textbf{Use of generative AI} During the preparation of this manuscript, the author used ChatGPT with the GPT-6 model (OpenAI) for English-language editing, typographical checks, assistance in checking calculations and logical consistency, and as a supplementary aid in some parts of the mathematical reasoning. All mathematical content, including any ideas or suggestions arising from these interactions, was independently examined and verified by the author. The author takes full responsibility for the content of the manuscript.

\end{document}